\documentclass[12pt]{amsart} \usepackage{latexsym, amssymb}
\usepackage[T1]{fontenc}
\usepackage[latin1]{inputenc}
\usepackage{pstricks}

\newtheoremstyle{slthm}
  {9pt}
  {5pt}
  {\slshape}
  {}
  {\bfseries}
  {.}
  {.5em}
  {\thmname{#1}\thmnumber{ #2}\thmnote{ (#3)}}

\newtheoremstyle{prcl}
  {9pt}
  {5pt}
  {\slshape}
  {}
  {\bfseries}
  {.}
  {.5em}
  {\thmname{#3}\thmnumber{ #2}}

\newtheoremstyle{prblm}
  {9pt}
  {5pt}
  {\rm}
  {}
  {\bfseries}
  {.}
  {.5em}
  {\thmname{#3}\thmnumber{ #2}}

\theoremstyle{slthm}
\newtheorem{thm}{Theorem}[section]
\newtheorem{lemma}[thm]{Lemma}
\newtheorem{prop}[thm]{Proposition}
\newtheorem{cor}[thm]{Corollary}

\theoremstyle{definition}

\newtheorem{df}[thm]{Definition}

\newtheorem{nrmk}[thm]{Remark}
\newtheorem{nrmks}[thm]{Remarks}
\newtheorem{expl}[thm]{Example}

\theoremstyle{remark}
\newtheorem*{rmk}{Remark}

\theoremstyle{prcl}
\newtheorem*{prclaim}{Proclaim}

\theoremstyle{prblm}

\newenvironment{renumerate}
        {
         \begin{enumerate}}
        {\end{enumerate}}

\newcounter{flexnummark}

\DeclareMathOperator{\cl}{cl}

\DeclareMathOperator{\ir}{int}

\DeclareMathOperator{\im}{Im}
\DeclareMathOperator{\re}{Re}
\DeclareMathOperator{\supp}{supp}

\DeclareMathOperator\Arg{Arg}

\newcommand{\mdots}{\dots}

\newcommand{\rest}[1]{|_{#1}}

\newcommand{\into}{\longrightarrow}

\renewcommand{\hat}{\widehat}
\renewcommand{\tilde}{\widetilde}
\renewcommand{\bar}{\overline}
\newcommand{\mdot}{\ \!\cdot\ \!}

\def\Ind#1#2{#1\setbox0=\hbox{$#1x$}\kern\wd0\hbox to 0pt{\hss$#1\mid$\hss}
\lower.9\ht0\hbox to 0pt{\hss$#1\smile$\hss}\kern\wd0}

\def\Notind#1#2{#1\setbox0=\hbox{$#1x$}\kern\wd0\hbox to 0pt{\mathchardef
\nn=12854\hss$#1\nn$\kern1.4\wd0\hss}\hbox to
0pt{\hss$#1\mid$\hss}\lower.9\ht0 \hbox to
0pt{\hss$#1\smile$\hss}\kern\wd0}

\newcommand{\set}[1]{\left\{#1\right\}}

\newcommand{\NN}{\mathbb{N}}

\newcommand{\RR}{\mathbb{R}}
\newcommand{\CC}{\mathbb{C}}

\newcommand{\curly}[1]{\mathcal{#1}}
\newcommand{\A}{\curly{A}}
\newcommand{\B}{\curly{B}}
\newcommand{\C}{\curly{C}}

\newcommand{\E}{\curly{E}}

\newcommand{\G}{\curly{G}}

\newcommand{\J}{\curly{J}}

\renewcommand{\O}{\curly{O}}

\newcommand{\Q}{\curly{Q}}
\newcommand{\R}{\curly{R}}

\renewcommand{\a}{\mathbf{a}}
\renewcommand{\b}{\mathbf{b}}
\newcommand{\f}{\mathbf{f}}
\newcommand{\g}{\mathbf{g}}
\renewcommand{\k}{\mathbf k}
\renewcommand{\l}{\mathbf l}
\newcommand{\m}{\mathbf{m}}
\newcommand{\p}{\mathbf{p}}
\newcommand{\q}{\mathbf{q}}
\renewcommand{\r}{\mathbf{r}}
\newcommand{\s}{\mathbf{s}}
\renewcommand{\t}{\mathbf{t}}
\newcommand{\frB}{\mathbf{B}}
\newcommand{\frL}{\mathbf{L}}

\newcommand{\Rtilde}{\,\tilde\RR\,}

\newcommand{\Rg}{\RR_{\G}}
\newcommand{\Rq}{\RR_{\Q}}

\newcommand{\Ps}[2]{\mathbb{#1}[\![#2]\!]}
\newcommand{\Pc}[2]{\mathbb{#1}\{#2\}}

\newcommand{\As}[2]{#1[\![#2]\!]}

\numberwithin{equation}{section}

\allowdisplaybreaks[2]

\begin{document}

\title{Transition maps at non-resonant hyperbolic
  singularities are o-minimal}

\author{T. Kaiser, J.-P. Rolin and P. Speissegger}

\address{Universit\"at Regensburg \\ NWF I-Mathematik \\ 93040
  Regensburg, \linebreak Germany}
\email{tobias.kaiser@mathematik.uni-regensburg.de}

\address{Universit\'e de Bourgogne \\ UFR Sciences et Techniques \\ 9
  avenue Alain Savary - B.P. 47870 \\ 21078 Dijon Cedex \\ France}
\email{Jean-Philippe.Rolin@u-bourgogne.fr} 

\address{McMaster University \\ Department of Mathematics \&
  Statistics, \linebreak 1280 Main Street West \\ Hamilton, Ontario L8S
  4K1 \\ Canada} \email{speisseg@math.mcmaster.ca}

\subjclass {37C27, 37E35, 03C64}

\keywords {Vector fields, transition maps, o-minimal structures}

\thanks {Supported by DFG, CNRS and NSERC}

\date{\today; Preprint.}

\begin{abstract}
  We construct a model complete and o-minimal expansion $\Rq$ of the
  field of real numbers such that, for any planar analytic vector
  field $\xi$ and any isolated, non-resonant hyperbolic singularity
  $p$ of $\xi$, a transition map for $\xi$ at $p$ is definable in
  $\Rq$.  This expansion also defines all convergent generalized power
  series with natural support and is polynomially bounded.
\end{abstract}

\maketitle
\markboth{T. Kaiser, J.-P. Rolin and P. Speissegger}{Transition maps
  at non-resonant hyperbolic singularities}

\section*{Introduction}

One of the motivations for this paper is the following: let $\xi$ be a
(real) analytic vector field on $\RR^2$ such that $\xi^{-1}(0) =
\{p\}$ is an isolated singularity of $\xi$.  We assume here that the
flow of $\xi$ near $p$ is as pictured in Figure 1 below:
there are two trajectories of $\xi$ at $p$, one incoming to $p$,
called $\gamma^-$, and the other outgoing from $p$, called $\gamma^+$.
To describe the flow of $\xi$ near these trajectories, we fix two
small segments $\Lambda^-$ and $\Lambda^+$ transverse to $\xi$ and
equipped with analytic charts $x$ and $y$ such that $x = 0$ is the
intersection point of $\gamma^-$ with $\Lambda^-$ and $x>0$ to the
right of $\gamma^-$, and similarly $y=0$ is the intersection point of
$\gamma^+$ with $\Lambda^+$ and $y>0$ above $\gamma^+$.  Then for all
sufficiently small $x>0$, the trajectory of $\xi$ crossing $\Lambda^-$
in the point $x$ later crosses $\Lambda^+$ in the point $y = g(x)$ of
$\Lambda^+$.

\begin{figure}[htbp]
  \begin{center}
    \input{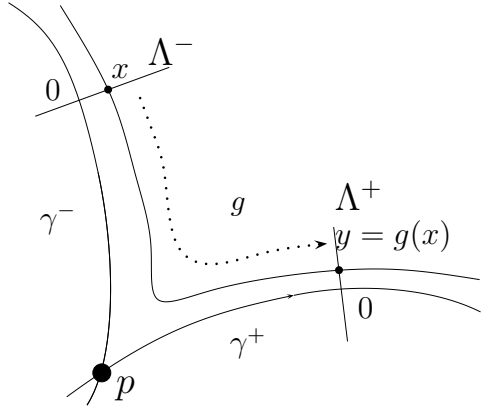}
  \end{center}
  \caption{\label{cap:Transition-map}Transition map $g$ at the
    hyperbolic singular point $p$}
\end{figure}

For any sufficiently small $\epsilon > 0$, the map $g:(0,\epsilon)
\into (0,\infty)$ defined in this way is called a \textbf{transition
  map} of $\xi$ at $p$.  The study of transition maps is at the heart
of Ilyashenko's solution of Dulac's Problem \cite{ily:dulac}.
Somewhat more precisely, Ilyashenko proves that any finite composition
of such transition maps has only finitely many isolated fixed points.
(Independently, Ecalle \cite{eca:dulac} proves that these maps are
\textit{analyzable} and deduces his own proof of Dulac's Problem.)

Ilyashenko's analysis of transition maps suggests to us the following
\medskip

\noindent\textbf{Question:} are the transition maps of $\xi$ near
$p$ definable in some \textit{fixed} \hbox{o-minimal} expansion $\R$
  of the real field?  \medskip

If the answer to this question is positive, it would follow that some
\textit{Poincar\'e return map near any polycycle of $\xi$} (see
Section \ref{family} for details) is also definable in $\R$, because
the family of functions definable in $\R$ is closed under composition.
It would then follow from Dulac's arguments \cite{dul:cycles} that
$\xi$ has at most finitely many limit cycles.

In this paper, we give a positive answer to the above question under
some restrictions on $\xi$.  First, we assume that the singularity $p$
of $\xi$ is \textbf{hyperbolic}, that is, the linear part of $\xi$ at
$p$ has two nonzero real eigenvalues of opposite signs.  In this
situation, Dulac proves in \cite[chapters 23 and 35]{dul:cycles} that
for a transition map $g$ as above, there exist a choice of charts $x$
and $y$, a $p_0>0$, real polynomials $p_j$ in one variable for $j = 1,
2, \dots$, and real numbers $0 < \nu_0 < \nu_1 < \cdots$, such that
$\lim_j \nu_j = +\infty$ and for every $n \in \NN$ the following
asymptotic relation holds:
\begin{itemize}
\item[(D)] $g(x) - p_0x^{\nu_0} - \sum_{j=1}^n p_j(\log x) x^{\nu_j} =
  o\left(x^{\nu_n}\right)$ as $x \to 0$.
\end{itemize}
Moreover, Ilyashenko obtains the following strengthening in Chapter 1
of \cite{ily:dulac}: a set $W \subseteq \CC$ is a \textbf{standard
  quadratic domain} if there are constant $c \in \RR$ and $C > 0$ such
that $$W = \set{z \in \CC:\ \re z < c - C \sqrt{|\im z|}}.$$ Then for
$g$ as above,
\begin{itemize}
\item[(I1)] there exists a standard quadratic domain $W \subseteq \CC$
  such that $g \circ \exp$ extends to a holomorphic mapping $G:W \into
  \CC$;
\item[(I2)] for all $n \in \NN$, we have the asymptotic relation
  \begin{equation*}
    G(z) - p_0 e^{\nu_0 z} - \sum_{j=1}^n p_j(z) e^{\nu_j z} =
    o\left(e^{\nu_n \re z}\right) \text{ as } |z| \to +\infty \text{
      in } W.
  \end{equation*}
\end{itemize}
Slightly abusing notations (to be clarified in Section
\ref{asymptotic_section} below), we summarize here conditions (I1) and
(I2) by saying that there is a standard quadratic domain $W \subseteq
\CC$ such that $g \sim_W p_0 x^{\nu_0} + \sum_{n=1}^\infty p_j(\log x)
x^{\nu_n}$.  A Phragmen-Lindel\"of argument \cite[p. 23]{ily:dulac}
shows that these conditions suffice to conclude that $g$ has at most
finitely many isolated fixed points.

The main body of Ilyashenko's proof consists in extending this
Phrag\-men-Lindel\"of argument to \textit{finite compositions} of
transition maps (not just in the hyperbolic case).  In contrast, our
approach is to try to prove that all transition maps generate an
o-minimal expansion of the real field.  Since finite compositions of
functions definable in an \hbox{o-minimal} structure are again
definable in that same structure, it would then follow that all finite
compositions of transition maps have finitely many isolated fixed
points.

In this paper, we carry out our approach under the additional
hypotheses of hyperbolicity and
\begin{itemize}
\item[(NR)] the singularity $p$ is \textbf{non-resonant}, that is, the
  ratio of the two eigenvalues of the linear part of $\xi$ at $p$ is
  an irrational number.
\end{itemize}
It follows from Dulac's argument that under the assumption (NR), the
polynomials $p_j$ in the asymptotic series of $g$ above are all
constant.  

Thus, we let $\Ps{R}{X^*}^\omega$ be the set of all formal power
series $F = \sum_{\alpha \ge 0} a_\alpha X^\alpha$ such that $a_\alpha
\in \RR$ for each $\alpha \ge 0$ and the support
\begin{equation*}
  \supp(F) := \set{\alpha \ge 0:\ a_\alpha \ne 0}
\end{equation*}
of $F$ is such that $\supp(F) \cap [0,R]$ is finite for every $R>0$.
Note that $\Ps{R}{X^*}^\omega$ is a subset of the set $\Ps{R}{X^*}$ of
all generalized power series defined by Van den Dries and Speissegger
in \cite{vdd-spe:genpower}.  (In the latter, o-minimality is
established for the expansion of the real field by all
\textit{convergent} generalized power series; in contrast, the
generalized power series studied here are in general not convergent.)
Since we do not use the larger class $\Ps{R}{X^*}$ here, we shall
routinely omit the superscript $\omega$.

Next, for every $\epsilon > 0$, we let $\Q_\epsilon$ be the set of all
functions $f:[0,\epsilon] \into \RR$ for which there exist an $F \in
\Ps{R}{X^*}$ and a standard quadratic domain $W \subseteq \CC$ such
that $(0,\epsilon] \subseteq \exp(W)$ and $f \sim_W F$.  (Thus by
definition, for any $f \in \Q_\epsilon$ the only point in
$[0,\epsilon]$ where $f$ is not necessarily analytic is $0$.)  Our
goal here is to prove that the expansion of the real field by all
functions in $\Q_1$ is o-minimal; to do so, we follow the method
developed in \cite{vdd-spe:genpower}.

Roughly speaking, the method in \cite{vdd-spe:genpower} goes as
follows: starting from a quasianalytic class $\C$ of functions in
several variables, we consider ``mixed'' functions that behave like
functions in $\C$ for some of the variables and are analytic in the
remaining variables \cite[section 5]{vdd-spe:genpower}.  We show that
the algebra of these mixed functions possesses certain closure
properties \cite[section 6]{vdd-spe:genpower}, most notably closure
under blow-up substitutions (which correspond to the charts of certain
blowings-up) and under Weierstrass preparation with respect to the
analytic variables.  These closure properties allow us \cite[section 7
and Proposition 8.4]{vdd-spe:genpower} to use resolution of
singularities to describe \textit{$\C$-sets}, which are defined by
equations and inequalities among functions from $\C$; in particular,
we show that $\C$-sets have finitely many connected components.  We
then adapt in \cite[section 8]{vdd-spe:genpower} Gabrielov's fiber
cutting argument to conclude that the complement of a projection of
a $\C$-set is again the projection of a $\C$-set.  The o-minimality
of the expansion of the real field by the functions in $\C$ follows as
outlined in \cite[section 2]{vdd-spe:genpower}.

This means in particular that we need to define classes $\Q_\rho$,
where $\rho \in (0,\infty)^m$ is a polyradius, of functions
$f:[0,\rho_1] \times \cdots \times [0,\rho_m] \into \RR$ with
analytic-extension and asymptotic properties in several variables
corresponding to (I1) and (I2) above.  It turns out, however, that the
most natural definition of these $\Q_\rho$ is insufficient to obtain
all the necessary closure properties, and we refer the reader to
Section \ref{division} for the correct definition of these classes.

Once the classes $\Q_\rho$ are introduced we define, for each $m \in
\NN$ and $m$-variable function $f \in \Q_{1, \dots, 1}$, a total
function $\tilde{f}:\RR^m \into \RR$ given by
\begin{equation*}
  \tilde{f}(x) := \begin{cases} f(x) &\text{if } x \in [0,1]^m, \\ 0
    &\text{otherwise.} \end{cases}
\end{equation*}
We let $\Rq := \left(\RR, <, 0, 1, +, -, \mdot, \left(\tilde{f}:\ f
    \in \Q_{1, \dots, 1}, m \in \NN\right)\right)$; note that in
particular, every function in $\Q_\epsilon$, for any $\epsilon > 0$,
is definable in $\Rq$.  Our partial answer to the question above is
the following:

\begin{prclaim}[Theorem A]
  The structure $\Rq$ is model complete, o-minimal and admits analytic
  cell decomposition. 
\end{prclaim}

Moreover, our method of proving Theorem A also gives a kind of
``Puiseux theorem'' for the one-variable functions definable in $\Rq$:

\begin{prclaim}[Theorem B]
  Let $\epsilon > 0$ and $f:(0,\epsilon) \into \RR$ be definable in
  $\Rq$.  Then there are a function $g \in \Q_\delta$ for some $\delta
  \in (0,\epsilon)$ and an $r \in \RR$ such that $g(0) \ne 0$ and
  $f(x) = x^r g(x)$ for all $x \in (0,\delta)$.
\end{prclaim}

Our discussion above of transition maps now implies:

\begin{prclaim}[Corollary]
  Assume that $p$ is a non-resonant hyperbolic singularity of $\xi$,
  and let $g:(0,\epsilon) \into (0,\infty)$ be a transition map of
  $\xi$ at $p$, expressed in the charts $x$ and $y$ such that (D)
  holds.  Then for every $\delta \in (0,\epsilon)$, the function
  $g\rest{(0,\delta)}$ is definable in $\Rq$.  \qed
\end{prclaim}

The transition maps discussed here are not the only functions of
interest definable in $\Rq$: in two forthcoming papers, the first
author shows that both the Riemann maps and the solutions of
Dirichlet's Problem on certain subanalytic domains with non-analytic
boundary are definable in $\Rq$.  

Moreover, we construct in Section \ref{family} of this paper an
\textit{analytic unfolding $\xi_\mu$} of $\xi$ in a neighbourhood of a
polycycle $\Gamma$ of $\xi$, in the spirit of Roussarie
\cite{rou:bifurcations}.  The unfolding $\xi_\mu$ is such that each
$\xi_\mu$ has the same set of singularities as $\xi$ and each of these
singularities is non-resonant hyperbolic.  Thus we obtain an analytic
family of transition maps at each of these singularities; in fact, we
show that each such family of transition maps is definable in $\Rq$.
It follows that there is a definable family of functions $P_\mu$ such
that for every $\mu$, the function $P_\mu$ is a Poincar\'e return map
for $\xi_\mu$ near $\Gamma$.  The uniform finiteness property of
o-minimal structures \cite[chapter 3, section 3]{vdd:vol1} then
implies that there is a uniform bound on the number of fixed points of
these functions $P_\mu$.  The question whether such uniform bounds
exist is related to Hilbert's 16th problem and remains open for
general $\xi$ and $\Gamma$.  (See \cite{ily:history} for a survey on
Hilbert's 16th problem and \cite{rou:bifurcations} for the
relationship between Hilbert's 16th problem and analytic unfoldings.)

Except for Section \ref{family}, the content of this paper is entirely
focussed on the construction of the o-minimal structure $\Rq$.  There
are various related questions that we do not know how to answer at
this point, such as:
\begin{enumerate}
\item Is $\Rg$, the o-minimal structure generated by all functions
  that are multisummable in the positive real direction
  \cite{vdd-spe:multisummable}, a reduct of $\Rq$?  Or is $\Rq$ a
  reduct of the Pfaffian closure of $\Rg$ (see \cite{spe:pfaffian})?
\item Are transition maps near \textit{resonant} hyperbolic
  singularities of $\xi$ definable in the Pfaffian closure of $\Rq$?
\end{enumerate}

This paper is organized as follows: in Section \ref{natural}, we
define the class of generalized power series with natural support, and
we establish some truncation properties needed later.  Our first and
most straightforward attempt at defining the classes $\Q_\rho$ is made
in Section \ref{asymptotic_section}, where we also establish some
useful criteria for functions to belong to these classes.  The reasons
why this first attempt is insufficient are explained in Section
\ref{tdte}, which lead us to a correct definition in Section
\ref{division}, where we also establish the first closure properties
needed to apply the method of \cite{vdd-spe:genpower}.  In fact, in
order to apply this method, we need to introduce corresponding
``mixed'' functions that are, roughly speaking, of asymptotic type
(I1) and (I2) in some variables and analytic in the other variables.
One of the reasons for studying these mixed classes is discussed
in Section \ref{weierstrass}, where we show that they admit
Weierstrass preparation with respect to the analytic variables.  To
reduce the study of sets defined by functions in the $\Q_\rho$ to that
of sets defined by functions in the mixed classes, we use the
blow-up substitutions introduced in Section \ref{blowups}.  Having
established all the properties necessary to apply the method of
\cite{vdd-spe:genpower}, we obtain Theorems A and B in Section
\ref{o-minimal}.

\section{Generalized power series with natural
  support}  \label{natural} 

Let $m \in \NN$, and let $X = (X_1, \dots, X_m)$ be a tuple of
indeterminates.  For $\alpha= (\alpha_1, \dots, \alpha_m) \in
[0,\infty)^m$, we write $X^\alpha := X_1^{\alpha_1} \cdots
X_m^{\alpha_m}$, and we let $X^*$ be the multiplicative monoid
consisting of all such $X^\alpha$, multiplied according to $X^\alpha
\cdot X^\beta = X^{\alpha+\beta}$.  The identity element of $X^*$ is
$X^0 = 1$, where $0 = (0, \dots, 0)$.  

Let $A$ be a commutative ring with $1 \ne 0$.  We let $\As{A}{X^*}$
denote the set of all formal power series $f(X) = \sum_{\alpha \ge 0}
a_\alpha X^\alpha$ such that $a_\alpha \in A$ for each $\alpha \ge 0$
and the \textbf{support} $\supp(f) := \set{\alpha \ge 0:\ a_\alpha \ne
  0}$ of $f$ is well-ordered, as defined in Section 4.4 of
\cite{vdd-spe:genpower}.  The elements of $\As{A}{X^*}$ are called
\textbf{generalized power series}.

\begin{df}
  \label{df:natural}
  A set $S \subseteq [0,\infty)^m$ is \textbf{natural} if
  $\Pi_{X_i}(S) \cap [0,R]$ is finite for every $R>0$ and each $i=1,
  \dots, m$, where $\Pi_{X_i}:\RR^m \into \RR$ is the projection on
  the coordinate $X_i$.

  We denote by $\As{A}{X^*}^\omega$ the set of all generalized power
  series in $X$ whose support is a natural subset of $[0,\infty)^m$.
\end{df}

\subsection*{Convention}
All results established in Section 4 of \cite{vdd-spe:genpower} go
through literally with $\As{A}{X^*}^\omega$ in place of $\As{A}{X^*}$,
as already pointed out in the second concluding remark of that paper.
To simplify notations, we shall from now on omit the superscript
$\omega$.  Thus, throughout this paper, every series in $\As{A}{X^*}$
is assumed to have natural support, and all results in
\cite{vdd-spe:genpower} referenced here are interpreted in this
context.  \medskip

Let $Y = (Y_1, \dots, Y_n)$ be another tuple of indeterminates.  For
$a,b \in \RR^k$, we write $a\le b$ (resp. $a < b$) if and only if $a_i
\le b_i$ (resp. $a_i < b_i$) for $i=1, \dots, k$.

\begin{df}
  \label{truncation_division}
  Let $S \subseteq [0,\infty)^{m} \times \NN^n$ and $F = \sum
  a_{(\alpha,\beta)} X^{\alpha} Y^\beta \in \As{A}{X^*,Y}$.  We define
  $\inf S = (a,b) = (a_1, \dots, a_m, b_1, \dots, b_n)$, where $a_i
  := \inf(\Pi_{X_i}(S))$ for $i=1, \dots, m$ and $b_i :=
  \min(\Pi_{Y_i}(S))$ for $i=1, \dots, n$, and we put
  \begin{equation*}
    F_{S} := \sum_{(\alpha,\beta) \in S}
    a_{(\alpha,\beta)} X^{\alpha-a} Y^{\beta-b}, 
    \quad\text{an element of } \As{A}{X^*,Y}.  
  \end{equation*}
  For $\gamma \in [0,\infty)^m \times \NN^n$, we write $F_\gamma$ in
  place of $F_{\{\alpha:\ \alpha \ge \gamma\}}$.
\end{df}

\begin{rmk}
  Let $F,G \in \As{A}{X^*,Y}$ and $\gamma \in [0,\infty)^m \times
  \NN^n$.  Then $(F+G)_\gamma = F_\gamma + G_\gamma$.
\end{rmk}

In the remainder of this section, we study how the operation $F
\mapsto F_\gamma$ behaves with respect to various other operations on
natural power series.  

\subsection*{Differentiation}

Let $F = \sum a_{(\alpha,\beta)} X^\alpha Y^\beta \in \As{A}{X^*,Y}$.
For $\beta \in \NN^n$ and $j=1, \dots, n$, we put $\beta^j :=
(\beta_1, \dots, \beta_{j-1}, \beta_j-1, \beta_{j+1}, \dots,
\beta_n)$.  We define
\begin{equation*}
  \partial_i F := \sum \alpha_i \cdot a_{(\alpha,\beta)} X^\alpha
  Y^\beta \quad\text{for } i=1, \dots, m,
\end{equation*}
and
\begin{equation*}
  \frac{\partial F}{\partial Y_j} := \sum \beta_j \cdot
  a_{(\alpha,\beta)} X^{\alpha} Y^{\beta^j}
  \quad\text{for } j = 1, \dots, n.
\end{equation*}
Note that each $\partial_i F$ and each $\partial F/\partial Y_j$
belongs to $\As{A}{X^*,Y}$.  Since each $\partial_i$ is a derivation
on $\As{A}{X^*,Y}$, we obtain:

\begin{lemma}
  \label{derivation_truncation}
  Let $F \in \As{A}{X^*,Y}$ and $\gamma \in [0,\infty)^m$.  Then:
  \begin{enumerate}
  \item $(\partial F/\partial Y_j)_{(\gamma,0)}
    = \partial (F_{(\gamma,0)})/\partial Y_j$ for each $j = 1, \dots, n$.
  \item $(\partial_i F)_{(\gamma,0)} = \gamma_i \cdot F_{(\gamma,0)} +
    \partial_i(F_{(\gamma,0)})$ for each $i=1, \dots, m$; in
    particular, if $\gamma_i = 0$, then $(\partial_i F)_{(\gamma,0)}
    = \partial_i(F_{(\gamma,0)})$.
  \end{enumerate}
\end{lemma}

\begin{proof}
  Parts (1) straightforward.  Since $X^\gamma \cdot G_\gamma$ is just
  the truncation of $G$ at $\gamma$ for any $G \in \As{A}{X^*}$, we have
  $X^\gamma \cdot (\partial_i F)_{(\gamma,0)} = \partial_i \left(X^\gamma
    \cdot F_{(\gamma,0)}\right)$.  Part (2) follows from the latter, because
  $\partial_i$ is a derivation.
\end{proof}

Before continuing with the behavior of the operation $F \mapsto
F_\gamma$, we make some crucial observations.

\subsection*{Representation}
Let $I \subset \{1,\dots,m\}$.  Below we write $X_I := (X_{i})_{i\in
  I}$, and if $x \in \RR^{m}$ we write $x_I:=(x_{i})_{i\in I}$.  We
let $\Pi_I:\RR^{m+n} \into \RR^{|I|}$ be the projection defined by
$\Pi_I(x,y) := x_I$.  We write $\bar{I} := \{1, \dots, m\} \setminus
I$.  For $F \in \As{A}{X^*,Y}$ and $\gamma \in [0,\infty)^m$, we put
$B_{\gamma,\emptyset} := \{0\}$, and for nonempty $I \subset
\{1,\dots,m\}$, we put $$B_{\gamma,I}=B_{\gamma,I}(F) :=
\set{\alpha\in \Pi_I(\supp F):\ \alpha < \gamma_I}.$$

\begin{lemma}
  \label{lemma:representation}
  Let $F \in \As{A}{X^*,Y}$ and $\gamma \in [0,\infty)^m$.  Then for
  each $I \subseteq \{1, \dots, m\}$, the set $B_{\gamma,I}$ is
  finite, and for each $\alpha \in B_{\gamma,I}$, there is a unique
  $F_{\gamma, I, \alpha} \in \As{A}{X_{\bar I}^*,Y}$ such that
  \begin{equation}
    \label{eq:representation}
    F(X,Y) = \sum_{I \subset \{1,\dots,m\}} X_{\bar I}^{\gamma_{\bar I}}
    \left(\sum_{\alpha \in B_{\gamma,I}} X_I^{\alpha}
      F_{\gamma,I,\alpha}(X_{\bar I},Y)\right).
  \end{equation}
\end{lemma}

\begin{proof}
  For each $(\alpha,\beta) \in \supp F$, there is a unique $I
  \subseteq \{1, \dots, m\}$ such that for all $i=1, \dots, m$, we
  have $\alpha_i < \gamma_i$ iff $i \in I$.  Moreover, $\alpha_I \in
  B_{\gamma,I}$ and $X^\alpha / \left(X_{\bar I}^{\gamma_{\bar I}}
    X_I^{\alpha_I}\right)$ is a monomial in $X_{\bar I}$.
\end{proof}

\begin{df}
  \label{df:representation}
  Let $F \in \As{A}{X^*,Y}$ and $\gamma \in [0,\infty)^m$.  We call
  the right-hand side of equation \eqref{eq:representation} the
  \textbf{$\gamma$-representation} of $F$.  For each $I \subseteq \{1,
  \dots, m\}$, we put
  \begin{equation*}
    F_{\gamma,I}(X,Y) := \sum_{\alpha \in B_{\gamma,I}} X_I^{\alpha}
    F_{\gamma,I,\alpha}(X_{\bar I},Y),
  \end{equation*}
  so that $F(X,Y) = \sum_{I \subseteq \{1, \dots, m\}} X_{\bar
    I}^{\gamma_{\bar I}} F_{\gamma,I}(X,Y)$.
\end{df}

\begin{nrmk}
  \label{rmk:elementary}
  In the situation of Definition \ref{df:representation}, we have
  $F_{(\gamma,0)} = F_{\gamma,\emptyset,0}$.  More generally, for each
  $\gamma \in [0,\infty)^m$, each $I \subseteq \{1, \dots, m\}$ and
  each $\alpha \in B_{\gamma,I}$, the series
  $F_{\gamma,I,\alpha}(X_{\bar I},Y)$, considered as an element of
  $\As{A}{X^*,Y}$, is equal to $F_S$ for some $S \subseteq
  [0,\infty)^m \times \NN^n$ of the following form:
\end{nrmk}

\begin{df}
  \label{elementary}
  Recall that a box in $\RR^k$ is a set of the form $B = \set{\alpha
    \in \RR^k:\ a_i \ast_{i,1} \alpha_i \ast_{i,2} b_i \text{ for each
    } i}$, where $a,b \in (\RR \cup \{\infty\})^k$ and $\ast_{i,1},
  \ast_{i,2} \in \{<,\le\}$.  A set $E \subseteq \RR^k$ is
  \textbf{elementary} if $E$ is a finite Boolean combination of boxes.
\end{df}

\begin{nrmk}
  \label{boolean}
  If $B \subseteq [0,\infty)^k$ is a box, then $[0,\infty)^k \setminus
  B$ is a finite union of pairwise disjoint boxes; also, any finite
  intersection of boxes is a box.  Therefore, every elementary set is
  a finite union of pairwise disjoint boxes.
\end{nrmk}

\begin{expl}
  \label{elementary_example}
  \begin{enumerate}
  \item For each $i \in \{1, \dots, k\}$ and every $a>0$, the set
    $\{\alpha \in [0,\infty)^k:\alpha_i = a\}$ is elementary.
  \item For all $a_1, \dots, a_k,b \ge 0$ and every natural set $S
    \subseteq [0,\infty)^k$, there is an elementary set $E \subseteq
    \RR^k$ such that $$S \cap \set{\alpha \in
      [0,\infty)^k:\sum_{i=1}^k a_i \cdot \alpha_i \ge b} = S \cap
    E.$$
  \end{enumerate}
\end{expl}

Intersections of natural sets and boxes can be further
simplified:

\begin{lemma}
  \label{natural_elementary}
  Let $B \subseteq [0,\infty)^k$ be a box and $S \subseteq
  [0,\infty)^k$ be a natural set.  Then there exist $\gamma, \delta^1,
  \dots, \delta^k \in [0,\infty)^k$ such that $\gamma \le \delta^j$
  for each $j$ and $S \cap B = S \cap \{\alpha \ge \gamma\} \setminus
  \bigcup_{j=1}^k S \cap \{\alpha \ge \delta^j\}$.
\end{lemma}

\begin{proof}
  Say $B = \set{\alpha \in \RR^k:\ a_i \ast_{i,1} \alpha_i \ast_{i,2}
    b_i \text{ for each } i}$, where $a,b \in (\RR \cup \{\infty\})^k$
  and $\ast_{i,1}, \ast_{i,2} \in \{<,\le\}$.  For each $i \in \{1,
  \dots, k\}$, we define $\gamma_i := \min\set{r \in \Pi_{X_i}(S): a_i
    \ast_{i,1} r}$ and $\delta_i := \max\set{r \in \Pi_{X_i}(S): r
    \ast_{i,2} b_i}$.  Then $S \cap B = S \cap \set{\gamma \le \alpha
    \le \delta}$, so the lemma follows with $\delta^j$ defined by
  $\delta^j_i := \gamma_i$ if $i \ne j$ and $\delta^j_j :=
  \min\{r \in \Pi_{X_j}(S):\ r > \delta_j\}$. 
\end{proof}

We now return to the study of the operation $F \mapsto F_\gamma$.
Since the observations below are rather technical and not clearly
motivated at this point, the reader may want to skip the rest of this
section and come back to it later as needed (while reading Section
\ref{division}, say).

\subsection*{Blow-up substitutions}

Let $i,j \in \{1, \dots, m\}$ be such that $i \ne j$, and let $\rho >
0$ and $\lambda \ge 0$.  Using the binomial expansion
\begin{equation*}
  (\lambda + X_j)^\beta := \sum_{p \in \NN} \begin{pmatrix} \beta \\
    p \end{pmatrix} \lambda^{\beta-p} X_j^p \quad\text{if } \lambda >
  0 \text{ and } \beta \ge 0,
\end{equation*}
we let $\frB^{\rho,\lambda}_{ij}:\As{A}{X^*} \into \As{A}{X^*}$ be the
unique $A$-algebra homomorphism satisfying
\begin{equation*}
  \frB^{\rho,\lambda}_{ij}(X_k) =
  \begin{cases}
    X_k &\text{if } k \ne i, \\ X_j^\rho (\lambda + X_i) &\text{if } k
    = i.
  \end{cases}
\end{equation*}
The homomorphism $\frB^{\rho,\lambda}_{ij}$ is called a \textbf{blow-up
  substitution}; we call $\frB^{\rho,0}_{ij}$ \textbf{singular}, and if
$\lambda > 0$, we call $\frB^{\rho,\lambda}_{ij}$ \textbf{regular}.  We
shall often write $\frB^{\rho,\lambda}_{ij} F$ in place of
$\frB^{\rho,\lambda}_{ij}(F)$, for $F \in \As{A}{X^*}$.

\begin{nrmk}
  \label{singular_expl}
  The substitution $\frB^{\rho,0}_{ij}$ is the $A$-algebra
  homomorphism $s^\rho_{ij}$ defined in Section 4.13 of
  \cite{vdd-spe:genpower}. 
\end{nrmk}

\begin{prop}
  \label{blowup_truncation}
  Let $F \in \As{A}{X^*}$ and put $X' := (X_1, \dots, X_{m-1})$.  
  \begin{enumerate}
  \item For $\lambda > 0$, we have $\frB^{\rho,\lambda}_{m,m-1} F \in
    \As{A}{(X')^*,X_m}$.
  \item Let $\gamma \in [0,\infty)^m$, and assume that $\gamma_m = 0$
    if $\lambda > 0$.  Then there are $\alpha_1, \dots, \alpha_k \in
    [0,\infty)^m$ and elementary sets $E_1, \dots, E_k \subseteq
    [0,\infty)^m$ such that
    \begin{equation*}
      \left(\frB^{\rho,\lambda}_{m,m-1} F\right)_\gamma = \sum_{l=1}^k
      X^{\alpha_l} \cdot \frB^{\rho,\lambda}_{m,m-1} F_{E_l}. 
    \end{equation*}
  \end{enumerate}
\end{prop}

\begin{proof}
  Part (1) follows from the binomial expansion.  Next, for $k=1,
  \dots, m$, we put
  \begin{equation*}
    \delta_k :=
    \begin{cases}
      \gamma_k &\text{if } k \ne m, \\ \max\set{\gamma_m,
        \frac{\gamma_{m-1}}{\rho}} &\text{if } k=m,
    \end{cases}
  \end{equation*}
  and $\delta := (\delta_1, \dots, \delta_m)$.  Let
  \begin{equation*}
    F(X) = \sum_{I \subset \{1, \dots, m\}} X_{\bar I}^{\delta_{\bar
        I}} \left( \sum_{\alpha \in B_{\delta,I}} X_I^\alpha
      F_{\delta,I,\alpha}(X_{\bar I}) \right)
  \end{equation*}
  be the $\delta$-representation of $F$.  Put $B'_{\delta,I} :=
  B_{\delta,I}$ if $I = \emptyset$ or $I = \{m-1\}$,
  $B'_{\delta,\{m\}} := \set{\alpha \in B_{\delta,I}:\ \alpha_m \ge
    \gamma_m}$ and $$B'_{\delta,\{m-1,m\}} := \set{\alpha \in
    B_{\delta,I}:\ \alpha_m \ge \gamma_m,\,\alpha_{m-1} + \rho
    \alpha_m \ge \gamma_{m-1}}.$$ Then by the hypothesis on $\gamma$,
  \begin{multline*}
    X^\gamma \cdot \left(\frB^{\rho,\lambda}_{m,m-1}F(X)\right)_\gamma
    \\ = \sum_{I \subseteq \{m-1,m\}} \sum_{\alpha \in B'_{\delta,I}}
    \frB^{\rho,\lambda}_{m,m-1} \left(X_{\bar I}^{\delta_{\bar I}}
      X_I^\alpha\right) \cdot \frB^{\rho,\lambda}_{m,m-1}
    (F_{\delta,I,\alpha}(X_{\bar I})).
  \end{multline*}
  Since each term $\frB^{\rho,\lambda}_{m,m-1} \left(X_{\bar
      I}^{\delta_{\bar I}} X_I^\alpha\right)$ on the right-hand side
  is divisible by $X^\gamma$, part (2) follows.
\end{proof}

\subsection*{Composition}

Let $F \in \As{A}{X^*,Y}$.  For the next lemma, we also let $\q =
(\q_1, \dots, \q_n) \in \NN^n$ and put $\k := |\q|$.  We let $Z =
(Z_1, \dots, Z_\k)$ be a tuple of indeterminates and define
\begin{equation*}
  F_{\q}(X,Z) := F\left(X, Z_1 + \cdots +
    Z_{\q_1}, \dots, 
    Z_{\q_1 + \cdots + \q_{n-1}+1} + \cdots + Z_\k\right).  
\end{equation*}
Note that $F_\q \in \As{A}{X^*,Z}$.  

\begin{prop}
  \label{comp_representation}
  Let $G = (G_1, \dots, G_n) \in (\As{A}{X^*,Y})^n$ be such that $G(0)
  = 0$.  Then the series $F(X,G(X,Y))$ belongs to $\As{A}{X^*,Y}$, and
  for each $\gamma \in [0,\infty)^m$, there are
  \begin{renumerate}
  \item a $p \in \NN$ and a tuple $\q \in \NN^n$ and, with
    $\k:= |\q|$,
  \item elementary $E_1, \dots, E_p \subseteq
    [0,\infty)^m \times \NN^\k$ and $B_{i,j} \subseteq [0,\infty)^{m}
    \times \NN^n$ for each pair $(i,j)$ satisfying $i \in \{1, \dots,
    n\}$ and $j \in \{1, \dots, \q_i\}$,
  \end{renumerate}
  such that, with $G_B := ((G_1)_{B_{1,1}}, \dots,
  (G_1)_{B_{1,\q_1}}, (G_2)_{B_{2,1}}, \dots,
  (G_n)_{B_{n,\q_n}})$, we have $G_B(0) = 0$, each term
  $(X,G_B(X,Y))^{\inf E_j}$ is divisible by $X^\gamma$ in
  $\Ps{A}{X^*,Y}$ and $$F(X,G(X,Y))_{(\gamma,0)} = \sum_{j=1}^p \frac
  {(X,G_B(X,Y))^{\inf E_j}}{X^\gamma} \cdot
  (F_\q)_{E_j}(X,G_B(X,Y)).$$
\end{prop}

\begin{proof}
  It is standard to check that $F(Y,G(X,Y)) \in \As{A}{X^*,Y}$; we
  leave the details to the reader.  Let $\gamma \in [0,\infty)^m$, and
  for each $j \in \{1, \dots, n\}$, we let
  \begin{equation*}
    G_j = \sum_{I \subseteq
      \{1,\dots,m\}} X_{\bar
      I}^{\gamma_{\bar I}}  \left(\sum_{\alpha \in 
        B_{I}(G_j)} X_I^\alpha \cdot
      (G_j)_{I,\alpha}(X_{\bar I},Y)\right) 
  \end{equation*}
  be the $\gamma$-representation of $G_j$.  (We omit the subscript
  $\gamma$ in these notations for the duration of this proof.)  We let
  $\J$ be the set of all triples $(j,I,\alpha)$ such that $j \in \{1,
  \dots, n\}$, $I \subseteq \{1, \dots, m\}$ and $\alpha \in
  B_{I}(G_j)$, and we put $\kappa:= |\J|$ and fix a bijection
  $\sigma:\{1, \dots, \kappa\} \into \J$.  Below, we write
  $\sigma(\lambda) = (j(\lambda),I(\lambda),\alpha(\lambda))$ for
  $\lambda = 1, \dots, \kappa$; without loss of generality, we may
  assume that $j(\lambda) \le j(\lambda')$ whenever $\lambda \le
  \lambda'$.

  For each $\lambda = 1, \dots, \kappa$, we let $Z_\lambda$
  be a new indeterminate and put
  \begin{equation*}
    G_\lambda(X,Y) := X_{\bar
      I}^{\gamma_{\bar I}} X_I^\alpha \cdot
    (G_j)_{I,\alpha}(X_{\bar I},Y)
  \end{equation*}
  with $(j,I,\alpha) = (j(\lambda),I(\lambda),\alpha(\lambda))$.  We
  write $G_B := (G_1, \dots, G_\kappa)$; note that $G_B(0) = 0$.  We
  also put
  \begin{equation*}
    H(X,Z) := F\left(X, \sum_{j(\lambda)=1} Z_\lambda\,,
      \dots, \sum_{j(\lambda)=n} Z_\lambda\right), 
  \end{equation*}
  and we write $H(X,Z) = \sum c_{\mu,\nu} X^\mu Z^\nu$, where $\mu$
  ranges over $[0,\infty)^m$ and $\nu$ ranges over $\NN^\kappa$.  Note
  that $F(X,G(X,Y)) = H(X,G_B(X,Y))$ and $H = F_{\q}$ for some
  $\q \in \NN^n$ with $|\q| = \kappa$.

  For $Q \subseteq \{1, \dots, \kappa\}$, we put $I_Q :=
  \bigcap_{\lambda \in Q} I(\lambda)$.  Let $\mu \in [0,\infty)^{m}$
  and $\nu \in \NN^\kappa$, and put $Q(\nu) := \set{\lambda:\
    \nu_\lambda \ne 0}$.  Then
  \begin{multline}
    \label{term_truncation}
    X^\gamma \cdot
    \left(X^\mu G_B(X,Y)^\nu\right)_{(\gamma,0)} = \\
    \begin{cases}
      X^\mu G_B(X,Y)^\nu &\text{if } \mu_i + \sum_{\lambda=1}^\kappa
      \nu_\lambda\cdot
      \alpha(\lambda)_i \ge \gamma_i \text{ for each } i \in I_{Q(\nu)}, \\
      0 &\text{otherwise.}
    \end{cases}
  \end{multline}
  Therefore, for each $Q \subseteq \{1, \dots, \kappa\}$, we let $S_Q
  \subseteq [0,\infty)^{m} \times \NN^k$ be the set defined by
  \begin{multline*}
    S_Q := \Big\{(\mu,\nu) \in \supp(H):\ \nu_\lambda = 0 \text{ iff }
    \lambda \notin Q,\\ \quad \mu_i + \sum_{\lambda \in Q}
    \nu_\lambda\cdot \alpha(\lambda)_i \ge \gamma_i \text{ for each }
    i \in I_Q \Big\}.
  \end{multline*}
  $S_Q$ is in turn the disjoint union of the following sets: for each
  map $\eta:I_Q \into Q$, we define
  \begin{multline*}
    S_{Q, \eta} := \Big\{ (\mu,\nu) \in S_Q:\ \mu_i +
    \sum_{\lambda \in Q,\,\lambda < \eta(i)} \nu_\lambda \cdot \alpha(\lambda)_i
    < \gamma_i \text{ and } \\ \mu_i + \sum_{\lambda \in Q,\,\lambda \le
      \eta(i)} \nu_\lambda \cdot \alpha(\lambda)_i \ge \gamma_i \text{ for
      each } i \in I_Q\Big\}.
  \end{multline*}
  Then $S_Q = S_{Q,0} \cup \bigcup_{\eta} S_{Q,\eta}$, where $S_{Q,0}
  := \set{(\mu,\nu) \in S_Q:\ \mu_{I_Q} \ge \gamma_{I_Q}}$.  Moreover
  for each $\eta:I_Q \into Q$, we write $\max \eta := \max_{i \in I_Q}
  \eta(i)$ and $J_\eta := \{\lambda < \max\eta:\ \alpha(\lambda)_i >
  0$ for some $i \in I_Q\}$.  Then the set $C_{Q,\eta} :=
  \set{(\mu_{I_Q},\nu_{J_\eta}):\ (\mu,\nu) \in S_{Q,\eta}}$ is finite,
  and $S_{Q,\eta}$ is the disjoint union of the sets $S_{Q,\eta} \cap
  E_{Q,\eta,\beta}$ as $\beta$ ranges over $C_{Q,\eta}$, where
  \begin{multline*}
    E_{Q,\eta,\beta} := \Big\{ (\mu,\nu):\
    (\mu_{I_Q},\nu_{J_\eta}) = \beta, \\
    \nu_{\max\eta} \cdot \alpha(\max\eta)_i \ge \gamma_i - \mu_i -
    \sum_{\lambda \in Q,\,\lambda < \max\eta} \nu_\lambda \cdot
    \alpha(\lambda)_i \text{ for each } i \in I_Q\Big\}.
  \end{multline*}

  By Example \ref{elementary_example}, each set $E_{Q,\eta,\beta}$ is
  elementary.  It follows from \eqref{term_truncation} above that
  \begin{multline*}
    X^\gamma \cdot F(X,G(X,Y))_{(\gamma,0)} = \\
    \sum_{Q,\eta,\beta} (X,G_B(X,Y))^{\inf E_{Q,\eta,\beta}} \cdot
    H_{E_{Q,\eta,\beta}}(X,G_B(X,Y)),
  \end{multline*} 
  and it follows from the definition of each $E_{Q,\eta,\beta}$ that
  $X^\gamma$ divides the factor $(X,G_B(X,Y))^{\inf E_{Q,\eta,\beta}}$
  in $\Ps{A}{X^*,Y}$, as required.
\end{proof}

\section{Natural asymptotic expansions}  \label{asymptotic_section}

Throughout this paper, we denote by $\|z\|$ the Euclidean norm of $z
\in \CC^n$, and we put $|z|:= z_1 + \cdots + z_n$ for such $z$.
We let $\frL = \set{(r,\varphi):\ r>0,\, \varphi \in \RR}$ be the
Riemann surface of the logarithm.  We fix an arbitrary
$m\in\mathbb{N}$ and write $$x = (x_1, \dots, x_m) = ((r_1,\varphi_1),
\dots, (r_m, \varphi_m))$$ for elements of $\frL^{m}$.  For such $x$,
we put $\|x\| := \|(r_1, \dots, r_m)\|$, $\arg(x) := (\varphi_1,
\dots, \varphi_m)$ and $$\log_m x := (\log r_1+i\varphi_1, \dots, \log
r_m + i\varphi_m) \in \CC^m;$$ we omit the subscript $m$ whenever it
is clear from context.  Note that $\log:\frL^m \into \CC^m$ is an
analytic isomorphism.  Below, we let $z = (z_1, \dots, z_m)$ range
over $\CC^m$.

Recall that for an open set $U \subseteq \frL$, a function $f:U \into
\CC$ is holomorphic if the function $f\circ\log^{-1}: \log(U) \into
\CC$ is holomorphic.  The set of holomorphic functions on $U$ is
denoted by $\O(U)$.

\begin{df}
  \label{df:holomorphic}
  Let $U \subseteq \frL^m$ be open.  For $f \in \O(U)$ and $i \in \{1,
  \dots, m\}$, we define $\partial_i f:U \into \CC$ by $$\partial_i
  f(x) := \frac{\partial (f \circ \log^{-1})}{\partial z_i} (\log
  x).$$

  Note that $\partial_i f \in \O(U)$.
\end{df}

\begin{expl}
  \label{expl:powers}
  Let $\alpha = (\alpha_{1},\dots,\alpha_{m}) \in \CC^{m}$.  We put
  $x^\alpha := x_{1}^{\alpha_{1}} \cdots x_{m}^{\alpha_{m}}$, where
  $x_i^{\alpha_i} := \exp(\alpha_i \log(x_i))$ for each $i$. The
  function $(\cdot)^{\alpha} : \frL^{m} \longrightarrow \CC$ is
  holomorphic and $\partial_i(x^\alpha) = \alpha_i x^\alpha$ for each
  $i$.  
\end{expl}

For $R>0$, we write $B_\frL(R):= \set{x \in \frL:\ \|x\| < R}$.

\begin{df}
  \label{quadratic}
  Let $W \subseteq \frL$.  The set $W$ is a \textbf{standard quadratic
    domain} if there are constants $c,C>0$ such that
  \begin{equation*}
    W = \set{(r,\varphi) \in \frL:\ 0 < r < c \exp\left( -C \sqrt{|\varphi|}
      \right)}. 
  \end{equation*}
\end{df}

Below, we put $\frL_0 := \frL \cup \{0\}$, and we extend the topology
on $\frL$ to a topology on $\frL_0$ by taking the set $B_\frL(R) \cup
\{0\}$, for $R>0$, as a basis of open neighborhoods of $0$ in
$\frL_0$.  For a subset $W$ of $\frL^m$, we shall write $\cl_0(W)$ for
the topological closure of $W$ in $\frL_0^m$.

\begin{rmk}
  Note that if $W \subseteq \frL$ is a standard quadratic domain, then
  $W \cup \{0\}$ is not an open neighborhood of $0$ in $\frL_0$.
  In particular, $\ir(\cl_0(W)) = W$.
\end{rmk}

\begin{df}
  \label{quadratic2}
  A set $W \subseteq \frL$ is a \textbf{quadratic domain} if $W$
  contains a standard quadratic domain and $\ir(\cl_0(W)) = W$.

  Let $U \subseteq \frL^m$ and $k \le m$.  We say that $U$ is
  \textbf{$k$-quadratic} if 
  \begin{renumerate}
  \item there is a quadratic domain $W \subseteq \frL$ and an $R>0$
    such that $W^k \times B_\frL(R)^{m-k} \subseteq U$;
  \item $0 \notin \ir(\cl_0(\Pi_{X_i} U))$ for each $i=1, \dots, k$.
  \end{renumerate}
\end{df}

\begin{nrmks}
  \label{rmk:quadratic}
  Let $k \le m$.
  \begin{enumerate}
  \item Let $R>0$ and $W \subseteq \frL$ be a quadratic domain such
    that $W \subseteq B_\frL(R)$.  If $l \ge k$, then $W^l \times
    B_\frL(R)^{m-l} \subseteq W^k \times B_\frL(R)^{m-k}$; therefore, every
    $k$-quadratic domain in $\frL^m$ contains an $l$-quadratic domain.
  \item Let $U \subseteq \frL^m$ be a $k$-quadratic domain, and let
    $r \in (0,\infty)^m$.  Then the set $V := \log^{-1}(\log(U) - r)$
    is a $k$-quadratic domain.
  \end{enumerate}
\end{nrmks}

We now fix an $m$-quadratic domain $U \subseteq \frL^m$.

\begin{df}
  \label{asymptotic}
  Let $f \in \O(U)$ and $F = \sum a_\alpha X^\alpha \in \Ps{C}{X^*}$.  We say
  that $f$ has \textbf{asymptotic expansion} $F$ on $U$ and write $f
  \sim_U F$, if for each $\a>0$ there is an $m$-quadratic domain $U_{\a}
  \subseteq U$ such that
  \begin{equation*}
    f(x) - \sum_{|\alpha| \le \a} a_{\alpha} x^{\alpha} =
    o\left(\|x\|^{\a}\right) \quad\text{as } \|x\| \to 0 \text{ on }
    U_\a.  
  \end{equation*}
  Note that in this situation, $F$ is the unique series $G \in
  \Ps{C}{X^*}$ with $f \sim_U G$; we therefore also write $Tf := F$,
  and we put $f(0) := \lim_{\|x\| \to 0,\,x \in U_\a} f(x)$ for any
  $\a > 0$.  We let $\A(U)$ be the set of all $f \in \O(U)$ for which
  there is an $F \in \Ps{C}{X^*}$ such that $f \sim_U F$.
\end{df}

\begin{nrmk}
  \label{O-rmk}
  \begin{enumerate}
  \item If $f \in \A(U)$ and $n \ge m$, then we consider $f$ as an
    holomorphic function $f:U \times \frL^{n-m}$ in the obvious way;
    under this identification, we get $\A(U) \subseteq \A(U
    \times \frL^{n-m})$.
  \item Let $f \in \O(U)$ and $V \subseteq U$.  Then $f\rest{V} \in
    \A(V)$ if and only if $f \in \A(U)$.
  \item The set $\A(U)$ is a $\CC$-algebra, and the map $T:\A(U) \into
    \Ps{C}{X^*}$ given by $T(f) := Tf$ is a $\CC$-algebra homomorphism
    satisfying $f(0) = (Tf)(0)$.
  \end{enumerate}
\end{nrmk}

For $R>0$, we  put
\begin{equation*}
  (0,R)_\frL := \set{x \in \frL:\ \|x\| < R \text{ and } \arg(x) = 0}.
\end{equation*}

\begin{prop}
  \label{asymptotic_to_0}
  The map $T:\A(U) \into \Ps{C}{X^*}$ is injective.
\end{prop}

\begin{proof}
  Let $f \in \A(U)$, and assume that $f \sim_U 0$; it suffices to show
  that $f = 0$.  Let $c,C>0$ and $$W := \set{(r,\varphi) \in \frL:\
    0<r<c \exp\left(-C\sqrt{|\varphi|}\right)}$$ be such that $W^m
  \subseteq U$.  For $s=(s_{1},\dots,s_{m}) \in (0,1]^{m}$, we define
  $f_{s}: W \into \CC$ by $f_s(r,\varphi) := f((s_{1}r,\varphi),
  \dots, (s_{m}r,\varphi))$.  Then $f_s \in \O(W)$, and our hypothesis
  implies that $f_s \sim_W 0$.  This means that
  \begin{equation*}
    \|f_s \circ \log^{-1}(z)\| = o(e^{-n\re z}) \quad\text{as } \re
    z \to -\infty \text{ in } \log(W) 
  \end{equation*}
  for every $n \in \NN$; hence $f_s = 0$ by Theorem 2 on p. 23 of
  \cite{ily:dulac}.  In particular, $f((s_1c,0), \dots, (s_mc,0))
  = 0$ for all $s \in (0,1)^m$, that is, $f \rest{(0,c)_\frL^{m}} = 0$.
  Therefore, the holomorphic map $h := f \circ \log^{-1}$ vanishes on
  $(-\infty,\log c)^{m}$.  Since $\log(U) \subseteq \CC^{m}$ is
  connected, it follows that $h=0$ and hence that $f=0$.
\end{proof}

\begin{prop}
  \label{prop:bad_derivatives}
  Let $f \in \A(U)$ and $i \in \{1, \dots, m\}$.  Then $\partial_i f$
  belongs to $\A(U)$ and satisfies $T(\partial_i f) = \partial_i(Tf)$.
\end{prop}

\begin{proof}
  Let $\a>0$, and assume that $U$ is an $m$-quadratic domain and
  $\|f(x)\| = o(\|x\|^\a)$ as $\|x\| \to 0$ with $x \in U$.  By Remark
  \ref{O-rmk}(2) and Example \ref{expl:powers}, it suffices to find an
  $m$-quadratic domain $V \subseteq U$ such that $\|\partial_i f(x)\|
  = o(\|x\|^\a)$ as $\|x\| \to 0$ with $x \in V$.

  We claim that $V := \log^{-1}(\log(U)-(1, \dots, 1))$ works; by
  Remark \ref{rmk:quadratic}(2), $V$ is an $m$-quadratic domain
  contained in $U$.  To see the claim, for each $r>0$, we put
  \begin{equation*}
    M_r := \max\set{\|f(x)\|:\ x \in U \text{ and } \|x\| \le r}
  \end{equation*}
  and
  \begin{equation*}
    N_r := \max\set{\|\partial_i f(x)\|:\ x \in V \text{ and } \|x\| \le
      r}. 
  \end{equation*}
  By assumption, we have $M_r / r^\a \to 0$ as $r \to 0$; we need to
  show that $N_r / r^\a \to 0$ as $r \to 0$.  To do so, it suffices to
  show that $N_r \le M_{r \cdot e}$ for each $r>0$, where $e :=
  \exp(1)$.  By the Cauchy estimates, we have for all $x \in V$ with
  $\|x\| \le r$ that
  \begin{equation*}
    \|\partial_if(x)\| = \left\|\frac{\partial (f \circ \log^{-1})} {\partial
      z_i} (\log x)\right\| \le M_{r \cdot e},
  \end{equation*}
  because $M_{r \cdot e}$ is the maximum of all $\|(f \circ
  \log^{-1})(z)\|$ such that $z \in \log(U)$ and $\re z_j \le \log r
  - 1$ for each $j$.  This finishes the proof of the proposition.
\end{proof}

We now also fix a $k \le m$.

\begin{df}
  \label{hol_proj}
  We let $\A^m_k(U)$ be the set of all $f \in \A(U)$ such that $Tf \in
  \Ps{C}{X_{\{1, \dots, k\}}^*, X_{\{k+1, \dots, m\}}}$.  We also let
  $\pi^m_k:\frL_0^m \into \frL_0^k \times \CC^{m-k}$ be the map
  defined by $\pi^m_k(x) = (y_1, \dots, y_m)$, where
  \begin{equation*}
    y_i :=
    \begin{cases}
      x_i &\text{if } i \le k, \\ x_i^1 &\text{if } i>k \text{
        and } x_i \ne 0, \\ 0 &\text{otherwise.}
    \end{cases}
  \end{equation*}
\end{df}
We also denote by $\cl_0$ the topological closure in $\frL_0^k \times
\CC^{m-k}$.  As usual, we shall omit the superscript $m$ if clear from
context.

Finally, for each $i = 1, \dots, m$, we let $p_i:\frL^m \into \frL^m$
be the map defined by $p_i(x) := y$ with $y_j:= x_j$ if $j \ne i$,
$\|y_i\| := \|x_i\|$ and $\arg(y_i):= \arg(x_i)+2\pi$.

\begin{prop}
  \label{prop:characterization}
  Let $f \in \A(U)$.  The following are equivalent:
  \begin{enumerate}
  \item $f \in \A_k(U)$.
  \item There are a $k$-quadratic domain $V \subseteq \frL^m$ and a
    holomorphic $f^\sharp : \ir(\cl_0(\pi_k(V))) \into \CC$ such that
    $f(x) = f^\sharp(\pi_k(x))$ for all $x \in U \cap V$.
  \item For every $i = k+1, \dots, m$ and every $x \in U$ satisfying
    $p_i(x) \in U$, we have $f(x) = f(p_i(x))$.
  \end{enumerate}
\end{prop}

\begin{proof}
  (2) $\Rightarrow$ (3): straightforward, since $U$ is connected.
  \smallskip 

  \noindent (3) $\Rightarrow$ (2): without loss of generality, we may
  assume that $U = W^m$ for some quadratic domain $W \subseteq \frL$.
  Let $R>0$ be such that $B(0,R) \subseteq \ir(\cl(\pi^1_0(W)))$, and
  put $V := W^k \times B_\frL(R)^{m-k}$.  Then the assumption on $f$ and
  Riemann's theorem on removable singularities imply that there is a
  holomorphic $g : \ir(\cl_0(\pi_k(V))) \into \CC$ such that $f(x) =
  g(\pi_k(x))$ for all $x \in U \cap V$, which proves (2).  \smallskip

  \noindent (3) $\Rightarrow$ (1): Assume (3) and let $i \in \{k+1,
  \dots, m\}$.  Write $Tf(X) = \sum a_\alpha X^\alpha$, and let
  $\alpha \in [0,\infty)^m$ be such that $\alpha_i \notin \NN$; we
  need to show that $a_\alpha = 0$.  Permuting coordinates if
  necessary, we may assume that $i=m$.  Let $R>0$ be such that the set
  $V := (0,R)_\frL^m \cup p_m((0,R)_\frL^m)$ is contained in $U$.

  Put $\a:= |\alpha|>0$, and define $g:V \into \CC$ by $g(x) := f(x) -
  \sum_{|\beta| \le \a} a_\beta x^\beta$.  Since $f \in \A(U)$, we
  have $g(x) = o\left(\|x\|^\a\right)$ as $\|z\| \to 0$ in $V$; in
  particular, for $s \in (0,R)$ and any $t \in (0,1)^m$, we have
  \begin{equation}
    \label{eq:diff_asym}
    g(st_1, \dots, st_m) - g(p_m(st_1, \dots, st_m)) = o(s^\a)
    \quad\text{as } s \to 0.
  \end{equation}
  On the other hand, 
  \begin{equation*}
    g(st_1, \dots, st_m) - g(p_m(st_1, \dots, st_m)) =
    \sum_{\substack{|\beta| \le \a \\ \beta_m \notin \NN}} a_\beta
    \left(e^{2\pi i \cdot 
        \beta_m}-1\right) t^\beta s^{|\beta|}.
  \end{equation*}
  It follows from \eqref{eq:diff_asym} that $\sum_{\substack{|\beta| \le
      \a \\ \beta_m \notin \NN}} a_\beta \left(e^{2\pi i \cdot
      \beta_m}-1\right) t^\beta = 0$.  Since $t \in (0,1)^m$ was
  arbitrary, we obtain that $a_\beta = 0$ for every $\beta \in
  [0,\infty)^m$ satisfying $|\beta|\le \a$ and $\beta_m \notin \NN$;
  in particular, $a_\alpha = 0$.  \smallskip

  \noindent (1) $\Rightarrow$ (3): assume that $Tf \in \Ps{C}{X_{\{1,
      \dots, k\}}^*, X_{\{k+1, \dots, m\}}}$, and let $i \in \{k+1,
  \dots, m\}$.  Then there is a quadratic domain $V \subseteq U$ such
  that $p_i(V) \subseteq U$, and we define $g:V \into \CC$ by $g(x) :=
  f(x) - f(p_i(x))$.  Then $g \in \O(V)$, and since $Tf \in
  \Ps{C}{X_{\{1, \dots, k\}}^*, X_{\{k+1, \dots, m\}}}$, we have $g
  \sim_V 0$.  Thus $g = 0$ by Proposition \ref{asymptotic_to_0}, which
  proves (3).
\end{proof}

Assume that $U$ is a $k$-quadratic domain and let $f \in \A_k(U)$.
Let $V$ and $f^\sharp$ be as in Proposition
\ref{prop:characterization}(2); by analytic continuation, we may
assume that $U = V$.  We extend $f$ to a function on
$\pi_k^{-1}(\ir(\cl_0(\pi_k(U))))$ by putting $f(x) :=
f^\sharp(\pi_k(x))$.  For example, we let $W \subseteq \frL$ be a
quadratic domain and $R>0$ be such that $W^k \times B_\frL(R)^{m-k}
\subseteq U$.  Then the value $f(x',0)$ is well defined for all
$x' \in W^k \times B_\frL(R)^{m-k-1}$, and we have:

\begin{cor}
  \label{value_at_0}
  The function $g:W^k \times B_\frL(R)^{m-k-1} \into \CC$ defined by $g(x')
  := f(x',0)$ belongs to $\A_k^{m-1}(W^k \times B_\frL(R)^{m-k-1})$ and
  satisfies $Tg(X') = Tf(X',0)$. \qed
\end{cor}

Moreover, we let $i \in \{k+1, \dots, m\}$.  Then the partial
derivative $\partial f^\sharp/\partial z_i:\pi_k(U) \into \CC$ is
defined as usual; using this, we define the partial derivative
$\partial f/\partial x_i:U \into \CC$ by
\begin{equation*}
  \frac {\partial f}{\partial x_i} (x) := \frac {\partial
    f^\sharp}{\partial z_i}(\pi_k(x)).
\end{equation*}

Using the Cauchy estimates (similarly as in the proof of Proposition
\ref{prop:bad_derivatives}) and Proposition \ref{asymptotic_to_0}, we
obtain:

\begin{cor}
  \label{good_derivatives}
  Assume that $U$ is a $k$-quadratic domain, and let $f \in \A_k(U)$.
  Then for every $i=k+1, \dots, m$, the partial derivative $\partial
  f/\partial x_i$ belongs to $\A_k(U)$ and satisfies $T(\partial
  f/\partial x_i) = \partial (Tf)/\partial X_i$; in particular, $x_i
  \cdot (\partial f/\partial x_i)(x) = \partial_i f(x)$ for all $x \in
  U$. \qed
\end{cor}

Finally, we establish some criteria for membership in $\A_k(U)$: we
fix an $l \in \{k, \dots, m\}$.  For $x \in \frL^m$ we write $Y :=
(X_1, \dots, X_{l})$, $y := (x_1, \dots, x_{l})$, $Z:= (X_{l+1},
\dots, X_m)$ and $z := (x_{l+1}, \dots, x_m)$.  We assume that $U =
W^k \times B_\frL(R)^{m-k}$ for some quadratic domain $W \subseteq \frL$
and some $R>0$, and we let $f \in \A_k(U)$.  By Proposition
\ref{prop:characterization}, for each $x = (y,z) \in U$ we have a
convergent power series representation
\begin{equation*}
  f(x) = \sum_{p \in \NN^{m-l}} a_p(y) z^p.
\end{equation*}
Writing $Tf = \sum_{\substack{p \in \NN^{m-l} \\ \alpha \in [0,\infty)^k
    \times \NN^{l-k}}} a_{p,\alpha} \cdot Y^\alpha Z^p$, we put
\begin{equation*}
  A_p(Y) := \sum_{\alpha \in [0,\infty)^k \times \NN^{l-k}} a_{p,\alpha} \cdot
  Y^\alpha, \quad\text{for each } p \in \NN^{m-l}.
\end{equation*}
Shrinking $U$ if necessary, there is an $M>0$ such that $\|f(x)\| \le
M$ for all $x \in U$.  Finally, we put $U' := W^k \times B_\frL(R)^{l-k}$.
Since $a_p(y) = \frac1{p!} \cdot \partial^p f/\partial (z)^p(y,0)$ for
all $p \in \NN^{m-l}$ and all $y \in U'$, we obtain from the Cauchy
estimates and Corollaries \ref{value_at_0} and \ref{good_derivatives}: 

\begin{prop}
  \label{coefficients}
  For each $p \in \NN^{m-l}$, the function $a_p:U' \into \CC$ belongs
  to $\A^l_k(U')$ and satisfies $Ta_p = A_p$ and $\|a_p(y)\| \le
  M/R^{|p|}$ for all $y \in U'$.  \qed
\end{prop}

The following converse to Proposition \ref{coefficients} is our
principal test for membership in $\A_k(U)$:

\begin{prop}
  \label{O-membership}
  Let $S \subseteq [0,\infty)^k \times \NN^{l-k}$ be natural, let
  $A>0$, and for each $p \in \NN^{m-l}$ let $b_p \in \A^l_k(U')$ be
  such that
  \begin{renumerate}
  \item $\supp(Tb_p) \subseteq S$;
  \item $||b_p(y)\| \le A/R^{|p|}$ for all $y \in U'$.
  \end{renumerate}
  Then the function $g:U \into \CC$ defined by $g(x) := \sum_{p \in
    \NN^{m-l}} b_p(y) \cdot z^p$ belongs to $\A_k(U)$ and satisfies
  $Tf = \sum_{p \in \NN^{m-l}} Tb_p(Y) \cdot Z^p$.
\end{prop}

\begin{proof}
  It follows from the assumptions that $g$ is holomorphic on $U$; it
  remains to show that $g \sim_U \sum_{p \in \NN^{m-l}} Tb_p(Y) Z^p$.
  Let $\a>0$; for each $p \in \NN^{m-l}$, we write $Tb_p = \sum
  b_{p,\alpha} Y^\alpha$ and define $\epsilon_p: U' \into \CC$ by
  \begin{equation*}
    \epsilon_p(y) := b_p(y) - \sum_{|\alpha| \le \a-|p|} b_{p,\alpha}
    y^\alpha. 
  \end{equation*}
  After shrinking $W$ if necessary, there is a constant $C>0$ such
  that $|\epsilon_p(y)| \le C\|y\|^{\a-|p|}$ whenever $|p| \le \a$ and
  $y \in U'$.  For every $x \in U$, we have
  \begin{equation*}
    g(x) - \sum_{|(\alpha,p)| \le \a} b_{p,\alpha} y^\alpha z^p
    = \sum_{|p| \le \a} \epsilon_p(y) z^p + \sum_{|p| > \a}
    b_p(y) z^p.
  \end{equation*}
  By the above, $\sum_{|p| \le \a} \epsilon_p(y) z^p =
  o\left(\|x\|^\a\right)$ as $\|x\| \to 0$ in $U$; so it suffices to
  show, after shrinking $U$ again if necessary, that $\sum_{|p| > \a}
  b_p(y) z^p = o\left(\|x\|^\a\right)$ as $\|x\| \to 0$ in $U$.  Let
  $\rho:= \inf\set{|p|-\a:\ |p| > \a} > 0$.  By assumption (ii), we have
  for $x \in U$ satisfying $\|z\| \le \frac \rho 2$ that
  \begin{equation*}
    \left|\sum_{|p| > \a} b_p(y) z^p \right| \le \frac{A}{R^{\a+\rho}}
    \left(\sum_{|p| \ge \a + \rho} 2^{\a+\rho-|p|}\right) \cdot
    \|z\|^{\a+\rho} = o\left(\|x\|^\a\right) 
  \end{equation*}
  as $\|x\| \to 0$, as required.
\end{proof}

The following criterion for membership in $\A_k(U)$ will be
useful in Section \ref{family}:

\begin{cor}
  \label{uniform_criterion}
  Let $f:U \into \CC$ be holomorphic, and let $W \subseteq \frL$ be a
  quadratic domain and $R>0$ be such that $W^k \times B_\frL(R)^{m-k}
  \subseteq U$.  Let $S \subseteq [0,\infty)^k$ be a natural set, and
  assume that for each $\alpha \in S$, there is a holomorphic function
  $a_\alpha:B_\frL(R)^{m-k} \into \CC$ such that
  \begin{renumerate}
  \item for every $z \in B_\frL(R)^{m-k}$, the function $f_z:W^k \into \CC$
    defined by $f_z(y):= f(y,z)$ belongs to $\A_k(W^k)$ with $Tf_z =
    \sum_{\alpha \in S} a_\alpha(z) y^\alpha$;
  \item for each $\nu>0$, there are constants $K_\nu, \epsilon_\nu >
    0$ and a quadratic domain $W_\nu \subseteq W$ such that $\left\|
      f(y,z) - \sum_{|\alpha| \le \nu} a_\alpha(z) y^\alpha \right\|
    \le K_\nu \|y\|^{\nu+\epsilon_\nu}$ for all $(y,z) \in W_\nu^k
    \times B_\frL(R)^{m-k}$.
  \end{renumerate}
  Then $f \in \A_k\left(W^k \times B_\frL(R)^{m-k}\right)$, and if
  $Ta_\alpha = \sum a_{\alpha,p} z^p$ for each $\alpha \in S$, then
  $Tf = \sum a_{\alpha,p} y^\alpha z^p$.
\end{cor}

\begin{proof}
  By assumption (i) and Proposition \ref{O-membership}, it suffices to
  show that the function $f^{(l)} : W^k \into \CC$ defined by
  $f^{(l)}(y) := (\partial^l f/\partial y^l)(y,0)$ belongs to
  $\A_k(W^k)$ for each $l \in \NN^{m-k}$.  But for any such $l$, any
  $\nu>0$ and any $y \in W_\nu$, we get from the Cauchy estimates that
  \begin{align*}
    \left\| f^{(l)}(y) - \sum_{|\alpha| \le \nu} a_\alpha^{(l)}(0)
      y^\alpha \right\| &= \left\| \frac{\partial^l}{\partial y^l}\left(
        f(y,z) - \sum_{|\alpha| \le \nu} a_\alpha(z) y^\alpha \right)
      (y,0)\right\| \\ &\le \frac {K_\nu}{R^{|l|}} \|y\|^{\nu+\epsilon_\nu},
  \end{align*}
  as required.
\end{proof}

\section{Truncation-division, Taylor expansion and composition in the holomorphic variables}
\label{tdte}

It will be convenient from now on to more explicitely separate the
holomorphic variables from the non-holomorphic ones.  Thus, we let
$m,n \in \NN$, and let $U \subseteq \frL^{m+n}$ be an $m$-quadratic
domain.  Below, we let $y = (y_1, \dots, y_n)$ range over $\frL^n$ and
$Y = (Y_1, \dots, Y_n)$ be indeterminates.

\begin{nrmk}
  \label{permutations}
  Let $\sigma$ be a permutation of $\{1, \dots, m\}$ and $\tau$ be a
  permutation of $\{1, \dots, n\}$.  We associate to $\sigma$ and
  $\tau$ the substitution automorphisms $\sigma,\tau:\Ps{C}{X^*,Y}
  \into \Ps{C}{X^*,Y}$ defined by $\sigma(X,Y) := (X_{\sigma(1)},
  \dots, X_{\sigma(m)},Y)$ and $\tau(X,Y):= (X,Y_{\tau(1)}, \dots,
  Y_{\tau(n)})$ and the maps $\sigma, \tau:\frL^{m+n} \into
  \frL^{m+n}$ defined by $\sigma(x,y) := (x_{\sigma(1)}, \dots,
  x_{\sigma(m)},y)$ and $\tau(x,y) := (x,y_{\tau(1)}, \dots,
  y_{\tau(n)})$. Then for every $f \in \A_m(U)$, we have $f \circ
  \sigma \in \A_m(\sigma^{-1}(U))$ with $T(f \circ \sigma) =
  \sigma(Tf)$ and $f \circ \tau \in \A_m(\tau^{-1}(U))$ with $T(f
  \circ \tau) = \tau(Tf)$.
\end{nrmk}

\subsection*{Truncation-division}  First, we show that the natural
operations of truncation and division by monomials in the $Y$
variables of $Tf$, where $f \in \A_m(U)$, lead to new functions in
$\A_m(U)$. 

\begin{prop}
  \label{holo_division}
  Let $f \in \A_m(U)$ and $\delta \in \NN^n$.  
  \begin{enumerate}
  \item Assume that $Tf = Y^\delta \cdot G$ for some $G \in
    \Ps{C}{X^*,Y}$.  Then there is a $g \in \A_m(U)$ such that $Tg =
    G$.
  \item There is an $f_{(0,\delta)} \in \A_m(U)$ such that
    $Tf_{(0,\delta)} = (Tf)_{(0,\delta)}$.
  \end{enumerate}
\end{prop}

\begin{proof}
  Below, we write $\hat{y}_j := (y_1, \dots, y_{j-1}, 0, y_{j+1},
  \dots, y_n)$ for each $j = 1, \dots, n$.  

  (1) Working by induction on $|\delta|$ and using Remark
  \ref{permutations}, we may assume that $\delta = (0, \dots, 0, 1)$.
  Given $h \in \A_m(U)$, the function $I(h):U \into \CC$ defined by
  \begin{equation*}
    I(h)(x,y) := \int_0^1 h(x,y_1, \dots, y_{n-1}, ty_n) dt
  \end{equation*}
  belongs to $\A_m(U)$, where $t \cdot (r,\varphi) := (tr,\varphi)$
  for $t>0$ and $(r,\varphi) \in \frL$.  Since $f(x,\hat{y}_n) = 0$,
  we get from the fundamental theorem of calculus that $f(x,y) = y_n
  \cdot I(\partial f/\partial y_n)(x,y)$ for all $(x,y) \in U$, so
  part (1) follows from Corollary \ref{good_derivatives}.

  (2) We define $h:U \into \CC$ by
  \begin{equation*}
    h(x,y) := f(x,y) - \sum_{j=1}^n \sum_{p=0}^{\delta_j-1} \frac1{p!}
    \frac {\partial^p f}{\partial y_j^p}(x,\hat{y}_j) \cdot y_j^p.
  \end{equation*}
  By Corollaries \ref{value_at_0} and \ref{good_derivatives}, the
  function $h$ belongs to $\A_m(U)$ and $Th = Y^\delta \cdot G$ for
  some $G \in \Ps{C}{X^*,Y}$.  Part (2) now follows from part (1).
\end{proof}

\begin{nrmk}
  \label{truncation_remark}
  We do not know whether, in the situation of Proposition
  \ref{holo_division}, a corresponding statements holds for all
  variables, that is, whether
  \begin{itemize}
  \item[$(\ast)_f$] for every elementary set $S \subseteq [0,\infty)^m
    \times \NN^n$, there are an $m$-quadratic $V \subseteq U$ and an
    $f_S \in \A_m(V)$ such that $T f_S = (Tf)_S$.
  \end{itemize}
  This is the first of two reasons to eventually restrict our
  attention to a subclass of $\A_m(U)$ introduced in Section
  \ref{division}.
\end{nrmk}

\subsection*{Taylor expansion}
Next, we establish a Taylor expansion result with respect to the $y$
variables.  Let $V \subseteq \frL^p$ be open; recall that a map $\f:V
\into \frL^p$ is holomorphic if the function $\log_p \circ \f \circ
\log_p^{-1}:\log_p(V) \into \CC^p$ is holomorphic.  Note that if $W
\subseteq \frL^q$ and $\g:W \into \frL^p$ is holomorphic such that
$\g(W) \subseteq V$, then the composition $\f \circ g:W \into \frL^p$
is holomorphic.

\begin{df}[Translation]
  \label{shift}
  Let $\lambda \in \CC$ be nonzero, and denote by $\Arg\lambda \in
  (-\frac\pi2, \frac\pi2]$ the standard argument of $\lambda$.  We
  put $$\CC_\lambda := \set{z \in \CC:\ \Arg\lambda - \frac\pi2 < \arg
    z < \Arg\lambda + \frac\pi2},$$ and we let $\l_\lambda:\CC_\lambda
  \into \frL$ be defined by $\l_\lambda(z) := (|z|,\arg z)$ and put
  $D(\lambda) := \l_\lambda(B(\lambda,|\lambda|))$.  We let
  $\t_\lambda:B_\frL(\lambda) \into D(\lambda)$ be the holomorphic map
  defined by $\t_\lambda(x) := \l_\lambda\left(\lambda +
    (x)^1\right)$.
  
  For completeness' sake, we also define $\t_0:\frL \into \frL$ by
  $\t_0(x) := x$.
\end{df}

For $\lambda \in \CC^p$ and $w \in \frL^p$, we write $\t_\lambda(w) :=
(\t_{\lambda_1}(w_1), \dots, \t_{\lambda_p}(w_p))$.  Abusing
notation, we identify the set $\big\{x \in \frL^p_0:\ -\pi < \arg
  x_i \le \pi$ for each $i \big\}$ with $\CC^p$.

\begin{nrmk}
  \label{translate_derivative}
  Let $\lambda \in \CC^{m+n} \cap \cl_0(U)$ and $f \in \O(U)$.  An
  elementary calculation shows that for each $i \in \{1, \dots, m\}$,
  we have
  \begin{equation*}
    \partial_i(f \circ \t_\lambda) = \Big((\partial_i f) \circ
    \t_\lambda\Big) \cdot \frac {x_i^1}{\lambda_i + x_i^1}.
  \end{equation*}
\end{nrmk}

We now let $l \in \{1, \dots, n\}$, and we write $y' := (y_1, \dots,
y_{n-l})$, \linebreak $Y' := (Y_1, \dots, Y_{n-l})$, $z = (z_1, \dots,
z_l):= (y_{n-l+1}, \dots, y_n)$ and $Z = (Z_1, \dots, Z_l) :=
(Y_{n-l+1}, \dots, Y_n)$. We assume that $U = W^m \times B_\frL(R)^n$ for
some quadratic domain $W \subseteq \frL$ and some $R>0$.  We let $f
\in \A_m(U)$ and $\lambda \in B(0,R)^l$, and we put $U':= W^m \times
B_\frL(R)^{n-l}$.

\begin{lemma}
  \label{real_constants}
  Assume $(\ast)_f$ holds.  Then the formal sum $Tf(X,Y',\lambda)$
  gives a series in $\Ps{C}{X^*,Y'}$.  Moreover, the function $(x,y')
  \mapsto f(x,y',\lambda):U' \into \CC$, denoted simply by
  $f(x,y',\lambda)$, satisfies $(\ast)_{f(x,y',\lambda)}$ with
  $$T(f(x,y',\lambda)) = (Tf)(X,Y',\lambda).$$
\end{lemma}

\begin{proof}
  We assume that $l=1$; the general case follows by induction on
  $l$.  Throughout this proof, we let $\gamma$ and $p$ range over
  $[0,\infty)^m \times \NN^{n-1}$ and $\NN$, respectively.  We write
  $f(x,y) = \sum a_p(x,y') z^p$ as a convergent power series in $z =
  y_n$ and $Tf = \sum a_{(\gamma,p)} (X,Y')^\gamma Z^p$.  From
  $(\ast)_f$ and Propositions \ref{coefficients} and
  \ref{asymptotic_to_0}, we see that $(\ast)_{a_p}$ holds (with $n-1$
  in place of $n$) for each $p \in \NN$.  Hence
  \begin{renumerate}
  \item $f_{(\gamma,0)}(x,y) = \sum (a_p)_\gamma(x,y') z^p$ for all
    $\gamma$ and all sufficiently small $(x,y) \in U$;
  \item $a_{(\gamma,p)} = (a_p)_\gamma(0,0)$ for all $\gamma$ and $p$.
  \end{renumerate}
  It follows from (ii) and Proposition \ref{coefficients} that
  $$Tf(X,Y',\lambda) = \sum_\gamma
  \left(\sum_p a_{(\gamma,p)} \lambda^p\right) (X,Y')^\gamma$$ belongs
  to $\Ps{C}{X^*,Y'}$, which proves the first assertion.

  Next, we let $\a>0$ and consider the finite set $$S_\a :=
  \set{\gamma \in \Pi_{m+n-1}(\supp Tf):\ |\gamma| \le \a}.$$ Note
  that, by (ii) above and Proposition \ref{coefficients}, for every
  $\gamma \in S_\a$ we have $f_{\{\gamma\} \times \NN}(x,y) = \sum_p
  a_{(\gamma,p)} z^p$ for all sufficiently small $(x,y) \in U$.
  Therefore, after shrinking $U$ if necessary, the function $$g(x,y)
  := f(x,y) - \sum_{\gamma \in S_\a} f_{\{\gamma\} \times \NN}(x,y)
  \cdot (x,y')^\gamma$$ belongs to $\A_m(U)$ and satisfies $Tg = Tf -
  \sum_{\gamma \in S_\a} (Tf)_{\{\gamma\} \times \NN} \cdot
  (X,Y')^\gamma$.  By the above, it follows that $$Tg(X,Y) = Tf(X,Y) -
  \sum_{\gamma \in S_\a} \sum_p a_{(\gamma,p)} \cdot (X,Y')^\gamma
  Z^p\ ;$$ in particular, $f(x,y',\lambda) - \sum_{|\gamma| \le \a}
  \left(\sum_p a_{(\gamma,p)} \lambda^p \right) \cdot (x,y')^\gamma =
  o(\|(x,y')\|^\a)$ as $\|(x,y')\| \to 0$ in $U'$.  Hence
  $f(x,y',\lambda) \in \A_m(U')$ with $T(f(x,y',\lambda)) =
  (Tf)(X,Y',\lambda)$.

  Finally, given an elementary set $S \subseteq [0,\infty)^m \times
  \NN^{n-1}$ and arguing as above with $f_{S \times \{0\}}$ in place
  of $f$, we see that $(\ast)_{f(x,y',\lambda)}$ holds.
\end{proof}

We put $R':= \min_{i=1, \dots, l} (R - |\lambda_i|)$ and $V := W^m
\times B_\frL(R)^{n-l} \times B_\frL(R')^l$.  From the Taylor expansion theorem
for holomorphic functions, Propositions \ref{coefficients} and
\ref{O-membership} and Lemma \ref{real_constants}, we obtain:

\begin{cor}
  \label{taylor}
  Assume $(\ast)_f$ holds.  Then the function $g:V \into \CC$ defined
  by $g(x,y) := f(x,y',\t_\lambda(z))$ satisfies $(\ast)_g$ with
  \begin{equation*}
    Tg(X,Y) = \sum_{p \in \NN^l} \frac1{p!} \frac
    {\partial^p (Tf)}{\partial Z^p}(X,Y',\lambda) \cdot
    Z^p.  \qed
  \end{equation*}
\end{cor}

In the situation of Corollary \ref{taylor}, we also write
$\t_{(0,\lambda)} f$ and $T_{(0,\lambda)} f$ for the function $g$ and
the series $Tg$, respectively.

\begin{nrmk}
  \label{taylor_remark}
  We are not aware of a statement corresponding to Corollary
  \ref{taylor}  for translation in the $x$ variables.  This is the
  second reason to restrict our attention to a subclass of $\A_m(U)$,
  done in Section \ref{division}.  
\end{nrmk}

\subsection*{Composition} Let $f \in \A_m(U)$.  Let $V \subseteq
\frL^{m+n}$ be an $m$-quadratic domain and let $g = (g_1, \dots, g_n)
\in \A_m(V)^n$.  Abusing notation, for $(x,y) \in V$ we write
$(x,g(x,y)) \in U$ to mean $(x, g(x,y)) \in \ir(\cl_0(\pi_m(U)))$, and
if the latter is the case, we also write $f(x,g(x,y))$ in place of
$f^\sharp(x,g(x,y))$.  

For the next lemma, we assume that $g(0,0) = 0$ and $(x, g(x,y))
\in U$ for all $(x,y) \in V$, and we define the holomorphic function
$h:V \into \CC$ by $ h(x,y) := f(x,g(x,y)).  $

\begin{prop}
  \label{A-composition}
  The function $h$ belongs to $\A_m(V)$ and $Th(X,Y) =
  Tf(X,Tg(X,Y))$.
\end{prop}

We will deduce this proposition from the following two special cases:

\begin{lemma}
  \label{comp1}
  Assume that $Tg(X,0) = Tg(X,Y)$.  Then $h \in \A_m(V)$ and $Th(X,Y)
  = Tf(X,Tg(X,Y))$.
\end{lemma}

\begin{proof}
  We write $Tf(X,Y) = \sum a_{\gamma,p} X^\gamma Y^p$ and $Tg_j(X) =
  \sum b_{j,\delta} X^\delta$; note that $b_{j,0} = 0$ for each $j$.
  Then $(Tf)(X,(Tg)(X,Y)) = \sum c_\alpha X^\alpha$, where for each
  $\alpha \in [0,\infty)^m$,
  \begin{multline*}
    \Sigma(\alpha) := \big\{(\gamma,p,\delta):\ \gamma \in
    \Pi_m(\supp(Tf)),\, p \in \NN^n, \text{ and } \\ \delta =
    (\delta^1, \dots, \delta^{|p|}) \in \left(\bigcup
      \supp(g_j)\right)^{|p|} \text{ with } \gamma + \delta^1 + \cdots
    \delta^{|p|} = \alpha\big\}
  \end{multline*}
  and
  \begin{equation*}
    c_\alpha := \sum_{(\gamma,p,\delta) \in \Sigma(\alpha)}
    a_{\gamma,p} \cdot \prod_{j=1}^n \prod_{l=p_1 + \cdots + p_{j-1} +
      1}^{p_1 + \cdots + p_l} b_{j,\delta^l}.
  \end{equation*}
  Since each $b_{j,0} = 0$, each set $\Sigma(\alpha)$ is finite; in
  fact, with
  \begin{equation*}
    q(r) := \sum_{i=1}^m \left| \Pi_{X_i}\left(\set{\beta \in
        \bigcup\supp(g_j):\ |\beta| \le r}\right) \right|,
  \end{equation*}
  we have $|p| \le q(|\alpha|)$ for all $(\gamma,p,\delta) \in
  \Sigma(\alpha)$. 

  Let now $\a > 1$, and for all suitable $(x,y) \in \frL^{m+n}$, we
  define
  \begin{equation*}
    f^\a(x,y) := f(x,y) - \sum_{|\gamma| + |p| \le \a + q(\a)}
    a_{\gamma,p} x^\gamma y^p
  \end{equation*}
  and
  \begin{equation*}
    g_j^\a(x) := g_j(x,0) - \sum_{|\delta| \le \a} b_{j,\delta}
    x^\delta, \quad\text{for } j=1, \dots, n.
  \end{equation*}
  Then $f^\a(x,y) = o\left(\|(x,y)\|^{\a+q(\a)}\right)$ as $\|(x,y)\| \to
  0$ in some $m$-quadratic domain, and $g_j^\a(x) = o(\|x\|^\a)$ for
  each $j$ as $\|x\| \to 0$ in some quadratic domain.  Thus,
  $f^\a(x,g^a(x)) = o(\|x\|^\a)$, and there is a polynomial $P(x) = \sum
  d_\beta X^\beta$ such that $|\beta| > \a$ whenever $d_\beta \ne 0$
  and
  \begin{equation*}
    P(x) = f(x,g(x,0)) - f^\a(x,g^\a(x)) - \sum_{|\alpha| \le \a}
    c_\alpha x^\alpha
  \end{equation*}
  for all sufficiently small $x$.  In particular, $f(x,g(x,0)) -
  \sum_{|\alpha| \le \a} c_\alpha x^\alpha = o(\|x\|^\a)$, which proves
  the lemma.
\end{proof}

\begin{lemma}
  \label{comp2}
  Assume that $Tg(X,0) = 0$.  Then $h \in \A_m(V)$ and $Th(X,Y)
  = Tf(X,Tg(X,Y))$.
\end{lemma}

\begin{proof}
  After shrinking $U$ and $V$ if necessary, we may assume that $f$ is
  bounded on $U$.  By Proposition \ref{coefficients}, we can write
  $f(x) = \sum a_p(x) y^p$ for all $(x,y) \in U$ and $g_j(x,y) = \sum
  b_{j,p}(x) y^p$ for all $(x,y) \in V$ and $j=1, \dots, n$, and there
  are an $m$-quadratic domain $W \subseteq \frL^m$ and constants
  $A,B>0$ such that $a_p, b_{j,p} \in \A_m(W^m)$ and $\|a_p(x)\|,
  \|b_{j,p}(x)\| \le AB^{|p|}$ for all $x \in W^m$, $p \in \NN^n$ and
  $j=1, \dots, n$.  Our assumption on $g$ implies that $b_{j,0} = 0$
  for each $j$.  Thus, after shrinking $W$ if necessary, there is an
  $R>0$ such that for all $(x,y) \in W \times B_\frL(R)^n$,
  \begin{equation*}
    h(x,y) = f(x,g(x,y)) = \sum_{p \in \NN^n} a_p(x) g(x,y)^p = \sum_{r
      \in \NN^n} c_r(x) y^r,
  \end{equation*}
  where, for $r \in \NN^n$ and $x \in W^m$, we put
  \begin{equation*}
    c_r(x) := \sum_{(p,q) \in \Sigma(r)} a_p(x) \cdot \prod_{i=1}^n
    \prod_{j=p_1 + \cdots + p_{i-1}+1}^{p_1 + \cdots + p_i} b_{i,q^j}(x)
  \end{equation*}
  with 
  \begin{multline*}
    \Sigma(r) := \big\{(p,q):\ p \in \NN^n \text{ with } |p| \le |r|,
    \text{ and } \\ q = (q^1, \dots, q^{|p|}) \in (\NN^n
    \setminus\{0\})^{|p|} \text{ with } q^1 + \cdots + q^{|p|} =
    r\big\}.
  \end{multline*}
  Note that each $\Sigma(r)$ is finite, because $|p| \le |r|$ for each
  $(p,q) \in \Sigma(r)$; we only need to consider such $p$, because
  each $b_{j,0} = 0$.)  Since $h$ is bounded and holomorphic on $W^m
  \times B_\frL(R)^n$, there are $C,D>0$ such that $\|c_r(x)\| \le
  CD^{|r|}$ for all $x \in W^m$ and $r \in \NN^n$.  Finally, writing
  $Tf(X,Y) = \sum_{p \in \NN^n} A_p(X) Y^p$ and $Tg_j(X,Y) = \sum_{q
    \in \NN^n} B_{j,p}(X) Y^q$ for $j=1, \dots, n$, it follows from
  Remark \ref{O-rmk}(3) that each $c_r$ belongs to $\A_m(W^m)$ and
  satisfies
  \begin{equation*}
    T c_r(X) = \sum_{(p,q) \in \Sigma(r)} A_p(X) \cdot \prod_{i=1}^n
    \prod_{j=p_1 + \cdots + p_{i-1}+1}^{p_1 + \cdots + p_i} B_{i,q^j}(x).   
  \end{equation*}
  The claim now follows from Proposition \ref{O-membership}, because
  $\supp(Tc_r) \subseteq \Pi_m\big(\supp Tf(X,Tg(X,Y))\big)$ for
  each $r \in \NN^n$ and the latter is a natural set by Proposition
  \ref{comp_representation}.
\end{proof}

\begin{proof}[Proof of Proposition \ref{A-composition}]
  We define $f'(x,z,y) := f(x,z+y)$, $g^0(x) := g(x,0)$ and $g'(x,y)
  := g(x,y) - g(x,0)$ for all suitable $x \in \frL^m$ and $y,z \in
  \frL^n$.  By Lemma \ref{comp2}, there is an $m$-quadratic $U'
  \subseteq \frL^{m+2n}$ such that $f' \in \A_m(U')$.  Note that
  $T(g^0)(X,0) = T(g^0)(X,Y)$ and $T(g')(X,0) = 0$.  Hence by Lemmas
  \ref{real_constants} and \ref{comp1}, there is an $m$-quadratic $V'
  \subseteq \frL^{m+n}$ such that $f'(x,g^0(x,y),y) \in \A_m(V')$, and
  by Lemma \ref{comp2}, there is an $m$-quadratic $V'' \subseteq
  \frL^{m+n}$ such that $f'(x,g^0(x,y),g'(x,y)) \in \A_m(V'')$.  Since
  $f(x,g(x,y)) = f'(x,g^0(x),g'(x,y))$ for all suitable $(x,y) \in
  \frL^{m+n}$, the proposition follows from Remark \ref{O-rmk}(2).
\end{proof}

\section{Blow-up substitutions in the non-holomorphic variables}
\label{blowups}

We continue to work with $m,n \in \NN$ and an $m$-quadratic domain $U
\subseteq \frL^{m+n}$.  For each real $\rho > 0$, the map
$\p^\rho:\frL \into \frL$ defined by
\begin{equation*}
  \p^\rho(r,\varphi) := (r^\rho,\rho\varphi)
\end{equation*}
is holomorphic, and the map $\m:\frL^2 \into \frL$ defined by
\begin{equation*}
  \m((r_1,\varphi_1),(r_2,\varphi_2)) := (r_1r_2, \varphi_1+\varphi_2)
\end{equation*}
is holomorphic.  Note that for all $x,x_1,x_2 \in \frL$, we have
$(\p^\rho(x))^1 = x^\rho$ for each $\rho > 0$ and $(\m(x_1,x_2))^1 =
(x_1)^1 \cdot (x_2)^1$.  

If $m \ge 2$ and $\rho \in (0,\infty)^m$, we define the holomorphic
map $\p^\rho:\frL^m \into \frL$ by induction on $m$:
\begin{equation*}
  \p^\rho(x) := \m\left(\p^{\rho'}(x'), \p^{\rho_m}(x_m)\right),
\end{equation*}
where $x' := (x_1, \dots, x_{m-1})$ and $\rho' := (\rho_1, \dots,
\rho_{m-1})$.

\begin{df}
  \label{singular_blowup}
  Let $m \ge 2$ and $i,j \in \{1, \dots, m\}$ be such that $i \ne j$,
  and let $\rho > 0$.  The \textbf{singular blowing-up} $\s^\rho_{ij}
  : \frL^{m+n} \into \frL^{m+n}$ is defined as $\s^\rho_{ij} (x,y) =
  (z,y)$, where
  \begin{equation*}
   z_k :=
   \begin{cases}
     x_k &\text{if } k \ne i, \\ \m(\p^\rho(x_j),x_i) &\text{if } k =
     i. 
   \end{cases}
  \end{equation*}
\end{df}

\begin{prop}
  \label{singular_prop}
  Let $ m \ge 2$ and $f \in \A_m(U)$.  Then the there is an
  $m$-quadratic $V \subseteq \frL^{m+n}$ such that $\s^\rho_{ij}(V)
  \subseteq U$ and the function $f \circ \s^\rho_{ij}$ belongs to
  $\A_m(V)$ and satisfies $T\left(f \circ \s^\rho_{ij}\right) =
  \frB^{\rho,0}_{ij}(Tf)$.
\end{prop}

\begin{proof}
  Without loss of generality, we may assume that $i = m$ and $j =
  m-1$.  Below, we write $\s$ and $\frB$ in place of $\s^\rho_{ij}$
  and $\frB^{\rho,0}_{ij}$.  Let $W \subseteq \frL$ be quadratic and
  $1 > R>0$ be such that $f \in \A_m(W^m \times B_\frL(R)^n)$; we may
  assume that
  \begin{equation*}
    W = \set{(r,\varphi) \in \frL:\ 0 < r <
      c\exp\left(-C\sqrt{|\varphi|}\right)} 
  \end{equation*}
  for some $c,C>0$ satisfying $c<R$.  We let $D :=
  C/\min\{\sqrt{\rho},1\}$ and put
  \begin{equation*}
    W' := \set{(r,\varphi) \in \frL:\ 0 < r <
      c\exp\left(-D\sqrt{|\varphi|}\right)} \subseteq W
  \end{equation*}
  and $V := (W')^m \times B_\frL(R)^n$; we claim that $\s(V) \subseteq U$.
  To see this, we write $x_k = (r_k,\varphi_k)$ for $k=1, \dots, m$.
  Then
  \begin{align*}
    \left\|\left(\s(x)\right)_m\right\| &\le c^{\rho+1}
    \exp\left( -D\left( \rho\sqrt{|\varphi_{m-1}|} +
        \sqrt{|\varphi_m|} \right)\right) \\ &\le c \exp\left(
      -D\left( \rho\sqrt{|\varphi_{m-1}|} + \sqrt{|\varphi_m|}
      \right)\right);
  \end{align*}
  since 
  \begin{align*}
    C \sqrt{|\rho\varphi_{m-1} + \varphi_m|} &\le C\left( \sqrt\rho
      \sqrt{|\varphi_{m-1}|} + \sqrt{|\varphi_m|}\right) \\ &\le D
    \left(\rho \sqrt{|\varphi_{m-1}|} + \sqrt{|\varphi_m|}\right),
  \end{align*}
  the claim follows.

 Since $f
  \circ \s$ is holomorphic on $V$, for each $\beta \in \NN^n$ the function
  $a_\beta:(W')^m \into \CC$ defined by
  \begin{equation*}
    a_\beta(x') := \frac1{\beta!} \left(\frac {\partial^\beta
        f}{\partial y^\beta} \circ \s\right) (x,0)  
  \end{equation*}
  is holomorphic.  By Proposition \ref{O-membership}, it suffices to
  show that $a_\beta \in \A((W')^m)$ for each $\beta$.  We
  fix $\beta \in \NN^n$ and write $T\left(\frac {\partial^\beta
      f}{\partial y^\beta}(X,0)\right) = \sum a_{\alpha} X^\alpha$,
  and we let $\a>0$.  Shrinking $W$ and $W'$ if necessary, we may
  assume by Lemma \ref{real_constants} that
  \begin{equation}
    \label{preblow_asymptotic}
    \left\|\frac {\partial^\beta f}{\partial y^\beta}(x,0) -
      \sum_{|\alpha| \le \a} a_{\alpha} 
      x^\alpha \right\|
    = o\left(\|x\|^\a\right) \quad\text{as } \|x\| \to 0
    \text{ in }  W^m.
  \end{equation}
  We now define $\rho:[0,\infty)^m \into [0,\infty)^m$
  by $$\rho(\alpha) := (\alpha_1, \dots, \alpha_{m-2},
  \alpha_{m-1}+\rho\alpha_m, \alpha_m).$$ Note that
  $\frB(T(\partial^\beta f/\partial Y^\beta))(X,0) = \sum a_{\alpha}
  X^{\rho(\alpha)}$ and $\s(x)^{\alpha} = x^{\rho(\alpha)}$ for all $x
  \in (W')^m$ and $\alpha \in [0,\infty)^m$.  Since $W \subseteq
  B_\frL(1)$, we have $\|\s(x)\| \le \|x\|$, so it follows from
  \eqref{preblow_asymptotic} that
  \begin{equation}
    \label{blow_asymptotic}
    \left\| \left(\frac {\partial^\beta
        f}{\partial y^\beta} \circ \s\right)(x,0) - \sum_{|\alpha|
        \le \a} a_{\alpha} x^{\rho(\alpha)} \right\|  =
    o\left(\|x\|^\a\right)
  \end{equation}
  as $\|x\| \to 0$ in $(W')^m$.  Finally, for $x \in W^m$ we have
  \begin{multline*}
    \left\|\left(\frac {\partial^\beta f}{\partial y^\beta} \circ
        \s\right) (x,0) - \sum_{|\rho(\alpha)| \le \a} a_{\alpha}
      x^{\rho(\alpha)} \right\| \\ \le \left\| \left(\frac
        {\partial^\beta f}{\partial y^\beta} \circ \s\right)(x,0) -
      \sum_{|\alpha| \le \a} a_{\alpha}
      x^{\rho(\alpha)} \right\|  +
    \left\|\sum_{|\alpha| \le \a < |\rho(\alpha)|} a_{\alpha}
      x^{\rho(\alpha)} \right\|.
  \end{multline*}
  The right-hand side above is $o\left(\|x\|^\a\right)$ as $\|x\| \to
  0$ in $(W')^m$, by \eqref{blow_asymptotic} and because the sum in
  the second summand is finite and each of its summands has an
  exponent $\gamma$ satisfying $|\gamma| > \a$.  This proves the
  proposition. 
\end{proof}

\begin{df}
  \label{regular_blowup}
  Let $m \ge 2$ and $\lambda > 0$.  The \textbf{regular blowing-up}
  $\r^{\rho,\lambda}: \frL^{m-1} \times B_\frL(\lambda) \times \frL^n \into
  \frL^{m+n}$ is defined as $\r^{\rho,\lambda}(x,y) = (z,y)$, where
  \begin{equation*}
    z_k :=
    \begin{cases}
      x_k &\text{if } k < m, \\ \m\left( \p^\rho(x_{m-1}),
        \t_\lambda(x_m) \right) &\text{if } k=m.
    \end{cases}
  \end{equation*}
\end{df}

\begin{prop}
  \label{regular_prop}
  Let $m \ge 2$, $\lambda > 0$ and $f \in \A_m(U)$.  Then there is an
  $(m-1)$-quadratic $V \subseteq \frL^{m-1} \times B_\frL(\lambda) \times
  \frL^n$ such that $\r^{\rho,\lambda}(V) \subseteq U$ and the
  function $f \circ \r^{\rho,\lambda}:V \into \CC$ belongs to
  $\A_{m-1}(V)$ and satisfies $T\left(f \circ \r^{\rho,\lambda}\right)
  = \frB^{\rho,\lambda}_{m,m-1}(Tf)$.
\end{prop}

\begin{proof}
  Below, we write $\r$ and $\frB$ in place of $\r^{\rho,\lambda}$ and
  $\frB^{\rho,\lambda}_{m,m-1}$.  Let $W' \subseteq \frL$ be quadratic
  and $\min\{1,\lambda\} > R>0$ be such that $U':= (W')^m \times
  B_\frL(R)^n \subseteq U$; we may assume that
  \begin{equation*}
    W' = \set{(r,\varphi) \in \frL:\ 0 < r <
      c\exp\left(-C\sqrt{|\varphi|}\right)} 
  \end{equation*}
  for some $c,C>0$ satisfying $c<R$.  We let 
  \begin{equation*}
    D:= \frac C{\min\{\sqrt\rho,1\}} \quad\text{and}\quad d :=
    \min\set{c, \left( 
        \frac c{2\lambda \exp\left(D \sqrt{\pi/2}\right)} \right)^{1/\rho}},
  \end{equation*}
  and we put
  \begin{equation*}
    W := \set{(r,\varphi) \in \frL:\ 0 < r <
      d\exp\left(-D\sqrt{|\varphi|}\right)} \subseteq W'
  \end{equation*}
  and $V := W^{m-1} \times B_\frL(R)^{n+1}$; we claim that $\r(V) \subseteq
  U$.  To see this, we write $x_k = (r_k,\varphi_k)$ for $k=1, \dots,
  m$.  Then
  \begin{align*}
    \left\|\left(\r(x)\right)_m\right\| &\le d^\rho \exp\left( -D
      \rho\sqrt{|\varphi_{m-1}|}\right) \cdot 2\lambda \\ &\le c
    \exp\left( -D\left( \rho\sqrt{|\varphi_{m-1}|} +
        \sqrt{\pi/2} \right)\right).
  \end{align*}
  Since $|\arg(\t_\lambda(w))| \le \pi/2$ for all $w \in \frL$, we
  also get
  \begin{align*}
    C \sqrt{|\rho\varphi_{m-1} + \arg(\t_\lambda(x_m))|} &\le C\left(
      \sqrt\rho \sqrt{|\varphi_{m-1}|} + \sqrt{\pi/2}\right) \\
    &\le D \left(\rho \sqrt{|\varphi_{m-1}|} +
      \sqrt{\pi/2}\right).
  \end{align*}
  The claim follows.

  We write $x' := (x_1, \dots, x_{m-1})$ for $x \in \frL^m$.  Since $f
  \circ \r$ is holomorphic on $V$, for each $p \in \NN$ and each
  $\beta \in \NN^n$ the function $a_{(p,\beta)}:(W')^{m-1} \into \CC$
  defined by
  \begin{equation*}
    a_{(p,\beta)}(x') := \frac1{p! \beta!} \frac {\partial^p((\partial
      f/\partial y^\beta) \circ \r)}{\partial
      x_m^p} (x',\lambda,0) 
  \end{equation*}
  is holomorphic.  Moreover, we put $X' := (X_1, \dots, X_{m-1})$, and
  $\alpha' := (\alpha_1, \dots, \alpha_{m-1})$, and we fix $p \in \NN$
  and $\beta \in \NN^n$ and write $T\left(\frac {\partial^\beta
      f}{\partial y^\beta}(X,0)\right) = \sum a_{\alpha} X^\alpha$.
  By the above and Proposition \ref{O-membership}, it now suffices to
  show that $a_{(p,\beta)} \in \A((W')^{m-1})$ with $Ta_{(p,\beta)} =
  \frac1{\beta!} A_{(p,\beta)}$,
  where $$A_{(p,\beta)}(X') := \sum_{\alpha} \begin{pmatrix} \alpha_m \\
    p \end{pmatrix} \lambda^{\alpha_m - p} a_{\alpha} (X')^{\alpha'}
  X_{m-1}^{\rho\alpha_m}.$$ Let $\a>0$, and choose $\a' >
  \a/\min\{\rho,1\}$.  Shrinking $W$ and $W'$ if necessary, we may
  assume by Lemma \ref{real_constants} that
  \begin{equation*}
    \left\|\frac {\partial^\beta f}{\partial y^\beta}(x,0) -
      \sum_{|\alpha| \le \a'} a_{\alpha} 
      x^\alpha \right\|
    = o\left(\|x\|^{\a'}\right) \quad\text{as } \|x\| \to 0
    \text{ in }  W^m.
  \end{equation*}
  Therefore,
  \begin{equation}
    \label{blow_asymptotic2}
    \left\| \left(\frac {\partial^\beta
          f}{\partial y^\beta} \circ \r\right)(x,0) -
      \sum_{|\alpha| \le \a'} a_{\alpha} \cdot
      (\r(x))^{\alpha} \right\|  =
    o\left(\|\r(x)\|^{\a'}\right)
  \end{equation}
  as $\|x\| \to 0$ in $(W')^m$.  We now define $\rho:[0,\infty)^m
  \into [0,\infty)^{m-1}$ by $\rho(\alpha) := (\alpha_1, \dots,
  \alpha_{m-2}, \alpha_{m-1}+\rho\alpha_m)$.  Note that formally
  $A_{(p,\beta)}(X') = \sum_{\alpha} \begin{pmatrix} \alpha_m \\
    p \end{pmatrix} \lambda^{\alpha_m - p} a_{\alpha}
  (X')^{\rho(\alpha)}$, and for all $x \in (W')^m$ and $\alpha
  \in [0,\infty)^m$ that $\r(x)^{\alpha} =
  (x')^{\rho(\alpha)} \cdot \sum_{q \in \NN} \begin{pmatrix} \alpha_m
    \\ q \end{pmatrix} \lambda^{\alpha_m-q} x_m^q$.  Differentiating
  \eqref{blow_asymptotic2}, it follows from the Cauchy estimates and
  our choice of $\a'$ that
  \begin{equation*}
    \left\| \beta! \cdot a_{(p,\beta)}(x') - \sum_{|\alpha| \le
        \a'} \begin{pmatrix} \alpha_m \\
        p \end{pmatrix} \lambda^{\alpha_m-p} a_{\alpha}
      (x')^{\rho(\alpha)} \right\|  =
    o\left(\|x'\|^\a\right)
  \end{equation*}
  as $\|x'\| \to 0$ in $(W')^{m-1}$.  Finally, for $x' \in (W')^{m-1}$ we
  have
  \begin{multline*}
    \left\| \beta! \cdot a_{(p,\beta)}(x') - \sum_{|\rho(\alpha)|
        \le
        \a} \begin{pmatrix} \alpha_m \\
        p \end{pmatrix} \lambda^{\alpha_m-p} a_{\alpha}
      (x')^{\rho(\alpha)} \right\| \\ \le \left\| \beta! \cdot
      a_{(p,\beta)}(x') - \sum_{|\alpha| \le
        \a'} \begin{pmatrix} \alpha_m \\
        p \end{pmatrix} \lambda^{\alpha_m-p} a_{\alpha}
      (x')^{\rho(\alpha)} \right\| \\ +
    \left\|\sum_{\substack{|\alpha| \le
          \a' \\ |\rho(\alpha)| > \a}} \begin{pmatrix}  \alpha_m \\
        p \end{pmatrix} \lambda^{\alpha_m-p} a_{\alpha}
      (x')^{\rho(\alpha)} \right\|.
  \end{multline*}
  The right-hand side above is $o\left(\|x'\|^\a\right)$ as $\|x'\|
  \to 0$ in $(W')^{m-1}$, by \eqref{blow_asymptotic2} and because the
  sum in the second summand is finite and each of its summands has an
  exponent $\beta$ satisfying $|\beta| > \a$, as $|\rho(\alpha)| \ge
  |\alpha| \min\{\rho,1\}$ for all $\alpha$.  This proves the
  proposition. 
\end{proof}

\section{The class $\Q$}
\label{division}

Let $m,n \in \NN$, and let $U \subseteq \frL^{m+n}$ be an
$m$-quadratic domain.  Below, we let $y = (y_1, \dots, y_n)$ range
over $\frL^n$ and $Y = (Y_1, \dots, Y_n)$ be indeterminates.

With Remarks \ref{truncation_remark} and \ref{taylor_remark} in mind,
we now restrict our attention to a subclass of $\A_m(U)$.  
Abusing notation, we identify $[0,\infty)$ with the set $\{0\} \cup
(0,\infty)_\frL \subseteq \frL_0$.

\begin{df}
  \label{Q}
  We define the class $\Q^{m+n}_m(U)$ to be the set of all $f \in
  \A_m(U)$ such that for every $\gamma \in [0,\infty)^m$, 
  \begin{itemize}
  \item[(TD)] there are an $m$-quadratic $V = V(f,\gamma)
    \subseteq U$ and an $f_{(\gamma,0)} \in \A_m(V)$ such that $T
    f_{(\gamma,0)} = (Tf)_{(\gamma,0)}$;
  \item[(TE)]for every $\kappa \in [0,\infty)^m$ with $(\kappa,0) \in
    \cl_0(U)$, there is an $m$-quadratic $W= W(f_{(\gamma,0)},\kappa)
    \subseteq \frL^{m+n}$ such that $(\t_\kappa(x),y) \in V$ for all
    $(x,y) \in W$ and the function $\t_{(\kappa,0)} f_{(\gamma,0)} : W
    \into \CC$ defined by $(\t_{(\kappa,0)} f_{(\gamma,0)})(x,y) :=
    f_{(\gamma,0)}(\t_\kappa(x),y)$ belongs to $\A_m(W)$.
  \end{itemize}
\end{df}

We shall omit the superscript $m+n$ whenever clear from context.

\begin{nrmks}
  \label{quasirmks}
  \begin{enumerate}
  \item By Proposition \ref{asymptotic_to_0}, for each $f \in \Q_m(U)$
    and each $\gamma \in [0,\infty)^m$, the function $f_{(\gamma,0)}$
    in the definition above is unique in $\A_m(V(f,\gamma))$.
  \item Let $f:U \into \CC$ be holomorphic, and let $V \subseteq U$ be
    an $m$-quadratic domain.  Then by Remark \ref{O-rmk}(2) and the
    above definition, $f\rest{V} \in \Q_m(V)$ iff $f \in \Q_m(U)$.
  \item By definition, the collection $\Q_1(U)$ is equal to the set of
    all $f \in \A_1(U)$ such that (TD) holds for every $\gamma \in
    [0,\infty)^m$; in particular, $\Q^1_1(U) = \A^1_1(U)$.  Moreover,
    $\A_0(U) = \Q_0(U)$ by Proposition \ref{holo_division} and
    Corollary \ref{taylor}.
  \item Let $\sigma$ be a permutation of $\{1, \dots, m\}$ and $\tau$
    be a permutation of $\{1, \dots, n\}$, and let $f \in \Q_m(U)$.
    Then $f \circ \sigma$ and $f \circ \tau$ belong to $\Q_m(U)$.
  \end{enumerate}
\end{nrmks}

For $p \in \NN$ and $q \in \{1, \dots, p\}$, we say that $\rho \in
[0,\infty)^p$ is \textbf{$q$-zero} if $\rho_1 = \cdots = \rho_q =
0$ and $\rho_{q+1}, \dots, \rho_p > 0$.  From Proposition
\ref{prop:characterization}, we obtain:

\begin{cor}
  \label{strongQ}
  Let $f \in \Q_m(U)$, and let $\gamma \in [0,\infty)^m$, $k \in \{1,
  \dots, m\}$, $\kappa \in [0,\infty)^m$ be $k$-zero such that
  $(\kappa,0) \in \cl_0(U)$ and $\sigma$ a permutation of $\{1, \dots,
  m\}$.  Then there is a $k$-quadratic domain $W \subseteq \frL^{m+n}$
  such that $\sigma(\t_\kappa(x),y) \in V(f,\gamma)$ for all $(x,y)
  \in W$ and the function $\t_{(\kappa,0)}(f_{(\gamma,0)} \circ
  \sigma)$ belongs to $\A_k(W)$. \qed
\end{cor}

\begin{df}
  \label{germs}
  Let $\E^{m+n}_m$ be the union of all $\Q_m(U)$ as $U$ ranges over
  the $m$-quadratic domains in $\frL^{m+n}$.  We define an equivalence
  relation $\equiv$ on $\E^{m+n}_m$ as follows: $f \equiv g$ if and
  only of there is an $m$-quadratic domain $U \subseteq \frL^{m+n}$
  such that $f \rest{U} = g\rest{U}$.  We let $\Q^{m+n}_m$ be the set
  of all $\equiv$-equivalence classes.
\end{df}

We shall omit the superscript $m+n$ whenever it is clear from context.
We will not distinguish between $f \in \Q_m(U)$ and its equivalence
class in $\Q_m$, which we also denote by $f$.  With this
identification, whenever $U \subseteq \frL^{m+n}$ is an $m$-quadratic
domain, we have $\Q_m(U) \subseteq \Q_m$.  Moreover, for every $f,g
\in \E_m$ such that $f \equiv g$, we have $Tf = Tg$; hence, the map $f
\mapsto Tf:\E_m \into \Ps{C}{X^*,Y}$ induces a map $f \mapsto Tf:\Q_m
\into \Ps{C}{X^*,Y}$.  Finally, for $r \ge 0$ we simply write $x^r$
for the germ of the function $x \mapsto x^r:\frL \into \CC$.

\begin{lemma}
  \label{Q-addition}
  \begin{enumerate}
  \item Let $f,g \in \Q_m$ and $a \in \CC$.  Then $f+g \in \Q_m$
    and $af \in \Q_m$.  
  \item If $m+n \ge l \ge m$, then $\Q_m \subseteq \Q_l$.
  \item Let $f \in \Q_m$ and $(\lambda,\mu) \in (0,\infty)^{m+n}$.
    Then the function $$f(\m(\lambda_1,x_1), \dots, \m(\lambda_m,x_m),
    \m(\mu_1,y_1), \dots, \m(\mu_n,y_n))$$ belongs to $\Q_m$.
  \end{enumerate}
\end{lemma}

\begin{proof}
  (1) Let $\gamma \in [0,\infty)^m$; then $(Tf)_{(\gamma,0)} +
  (Tg)_{(\gamma,0)} = (T(f+g))_{(\gamma,0)}$ and $a (Tf)_{(\gamma,0)}
  = T(af)_{(\gamma,0)}$.  Let also $\kappa \in [0,\infty)^m$ be
  sufficiently small.  Then $f \circ \t_{(\kappa,0)} + g \circ
  \t_{(\kappa,0)} = (f+g) \circ \t_{(\kappa,0)}$ and $a (f \circ
  \t_{(\kappa,0)}) = (af) \circ \t_{(\kappa,0)}$, so we define
  $\t_{(\kappa,0)}(f+g) := \t_{(\kappa,0)}f + \t_{(\kappa,0)}g$ and
  $\t_{(\kappa,0)}(af) := a \t_{(\kappa,0)} f$.

  (2) Let $m+n \ge l \ge m$, and let $f \in \Q_m$ and $\rho \in
  [0,\infty)^l$.  Let also $\rho'$ be the least $\tau \in [0,\infty)^m
  \times \NN^{l-m}$ such that $\tau \ge \rho$.  Then we can take
  $f_{(\rho,0)} := (x,y)^{\rho'-\rho} \cdot f_{(\rho',0)}$.  It
  follows easily that $f \in \Q_l$.

  (3)   Writing $\m((\lambda,\mu),(x,y)):= (\m(\lambda_1,x_1), \dots,
  \m(\mu_n,y_n))$, and writing $(\lambda,\mu) \cdot (X,Y) :=
  (\lambda_1 X_1, \dots, \mu_n Y_n)$, we see that
  \begin{equation*}
    Tf((\lambda,\mu) \cdot (X,Y))_{(\gamma,0)} =
    (Tf)_{(\gamma,0)}((\lambda,\mu) \cdot (X,Y))
  \end{equation*}
  for all $\gamma \in [0,\infty)^m$, and that
  \begin{equation*}
    \t_{(\kappa,0)}(f(\m((\lambda,\mu),(x,y)))) =
    (\t_{(\lambda\kappa,0)} f)(\m((\lambda,\mu),(x,y)))
  \end{equation*}
  for all sufficiently small $(x,y) \in \frL^{m+n}$, where
  $\lambda\kappa := (\lambda_1\kappa_1, \dots, \lambda_m\kappa_m)$.
  Part (3) follows.
\end{proof}
 
\begin{prop}
  \label{Q-basics}
  Let $f \in \Q_m$.  
  \begin{enumerate}
  \item For every elementary set $S \subseteq [0,\infty)^m \times
    \NN^n$, there is a unique $f_S \in \Q_m$ such that $T(f_S) =
    (Tf)_S$.
  \item For every $k \in \{1, \dots, m\}$, every sufficiently small
    $k$-zero $(\kappa,\lambda) \in [0,\infty)^m \times \RR^n$ and
    every permutation $\sigma$ of $\{1, \dots, m\}$, the function
    $\t^\sigma_{(\kappa,\lambda)} f := f \circ \sigma \circ
    \t_{(\kappa,\lambda)}$ belongs to $\Q_k$.
  \end{enumerate}
\end{prop}

\begin{proof}
  (1) By Remark \ref{boolean} and Lemmas \ref{natural_elementary} and
  \ref{Q-addition}, it suffices to consider $S = \set{((\alpha,\beta)
    \in [0,\infty)^m \times \NN^n:\ (\alpha,\beta) \ge
    (\gamma,\delta)}$ for some $(\gamma,\delta) \in [0,\infty)^m
  \times \NN^n$.  By Lemma \ref{Q-addition} and Proposition
  \ref{asymptotic_to_0}, we may even assume that either $\gamma = 0$
  or $\delta = 0$.  We assume first that $\delta = 0$ and let
  $f_{(\gamma,0)}$ be as in (TD); we need to show that $f_{(\gamma,0)}
  \in \Q_m$.  So let $\gamma' \in [0,\infty)^m$ and $\kappa \in
  [0,\infty)^m$ be sufficiently small.  Since
  $(Tf)_{(\gamma+\gamma',0)} = (Tf_{(\gamma,0)})_{(\gamma',0)}$, we
  can take $(f_{(\gamma,0)})_{(\gamma',0)} := f_{(\gamma+\gamma',0)}$
  and $\t_{(\kappa,0)}\big((f_{(\gamma,0)})_{(\gamma',0)}\big)
  := \t_{(\kappa,0)} f_{(\gamma+\gamma',0)}$.  Second, the case
  $\gamma = 0$ follows from Proposition \ref{holo_division} and
  Corollary \ref{taylor}.

  (2) Let $k \in \{1, \dots, m\}$ and $(\kappa,\lambda) \in
  [0,\infty)^m \times \RR^n$ be sufficiently small and $k$-zero; by
  Remark \ref{quasirmks}(4), it suffices to prove that
  $\t_{(\kappa,\lambda)} f$ belongs to $\Q_k$.  Since
  $\t_{(\kappa,\lambda)} f = \t_{(\kappa,0)}(\t_{(0,\lambda)} f)$, we
  may assume by Corollary \ref{taylor} that $\lambda = 0$.  Let
  $\gamma \in [0,\infty)^k \times \{0\}^{m-k}$.  By (1), there is for
  each $I \subseteq \{1, \dots, m\}$ and each $\alpha \in B_{\gamma,I}
  = B_{\gamma,I}(Tf)$ a unique $f_{\gamma,I, \alpha} \in \Q_m$ such
  that
  \begin{equation*}
    f = \sum_{I \subseteq \{1, \dots, m\}} x_{\bar
      I}^{\gamma_{\bar I}} \left( \sum_{\alpha \in B_{\gamma,I}}
      x_I^\alpha \cdot f_{\gamma,I,\alpha} \right)
  \end{equation*}
  and each $f_{\gamma,I,\alpha}$ depends only on the variables
  $x_{\bar I}$ and $y$.  Since $\gamma_{k+1} = \cdots = \gamma_m = 0$,
  we have $B_{\gamma,I} = \emptyset$ whenever $I \nsubseteq \{1,
  \dots, k\}$.  Therefore,
  \begin{equation*}
    \t_{(\kappa,0)} f = \sum_{I \subseteq \{1, \dots, m\}} x_{\bar
      I}^{\gamma_{\bar I}} \left( \sum_{\alpha \in B_{\gamma,I}}
      x_I^\alpha \cdot \t_{(\kappa,0)} f_{\gamma,I,\alpha} \right),
  \end{equation*}
  and hence $T(\t_{(\kappa,0)} f) = \sum_I X_{\bar I}^{\gamma_{\bar
      I}} \left( \sum_\alpha X_I^\alpha \cdot T(\t_{(\kappa,0)}
    f_{\gamma,I,\alpha}) \right)$ is the unique
  $\gamma$-representation of $T(\t_{(\kappa,0)} f)$.  Since
  $f_{\gamma,\emptyset,0} = f_{(\gamma,0)}$, it follows that we can
  take $(\t_{(\kappa,0)} f)_{(\gamma,0)} := \t_{(\kappa,0)}
  f_{(\gamma,0)}$.  Moreover, if $\kappa' \in [0,\infty)^k \times
  \{0\}^{m-k}$ is sufficiently small, then
  $$\t_{(\kappa',0)}((\t_{(\kappa,0)} f)_{(\gamma,0)}) =
  \t_{(\kappa',0)}(\t_{(\kappa,0)} f_{(\gamma,0)}) =
  \t_{(\kappa+\kappa',0)} f_{(\gamma,0)},$$ so part (2) follows.
\end{proof}

From Proposition \ref{Q-basics} and Lemma \ref{lemma:representation},
we obtain:

\begin{cor}
  \label{cor:Q-basics}
  Let $f \in \Q_m$ and $\gamma \in [0,\infty)^m$.  Then for each
  $I \subset \{1, \dots, m\}$ and each $\alpha \in B_{\gamma,I} =
  B_{\gamma,I}(Tf)$, there is a unique $f_{\gamma,I,\alpha} \in \Q_m$
  such that
    \begin{equation*}
      f = \sum_{I \subset \{1, \dots, m\}} x_{\bar I}^{\gamma_{\bar
          I}} \left( \sum_{\alpha \in B_{\gamma,I}} x_I^\alpha \cdot 
        f_{\gamma,I,\alpha} \right)
    \end{equation*}
    and each $f_{\gamma,I,\alpha}$ depends only on the variables
    $x_{\bar I}$ and $y$. \qed
\end{cor}

\begin{prop}
  \label{Q-derivatives}
  Let $f \in \Q_m$.
  \begin{enumerate}
  \item For each $i=1, \dots, m$, the function $\partial_i f$ belongs to
    $\Q_m$ and satisfies $T(\partial_i f) = \partial_i(Tf)$.
  \item For each $j = 1, \dots, n$, the function $\partial f/\partial
    y_j$ belongs to $\Q_m$ and satisfies $T(\partial f/\partial y_j)
    = \partial (Tf)/\partial Y_j$.
  \item The function $g := f(x,y_1, \dots, y_{n-1},0)$ belongs to
    $\Q^{m+n-1}_m$.
  \end{enumerate}
\end{prop}

\begin{proof}
  Let $\gamma \in [0,\infty)^m$.  It follows from Propositions
  \ref{prop:bad_derivatives} and \ref{Q-basics} and Lemmas
  \ref{derivation_truncation} and \ref{Q-addition} that $(\partial_i
  f)_{(\gamma,0)} := \gamma_i \cdot f_{(\gamma,0)}
  + \partial_i(f_{(\gamma,0)})$ belongs to $\A_m(V)$ for a suitable $V
  \subseteq \frL^{m+n}$ and satisfies $T((\partial_i f)_{(\gamma,0)})
  = \partial_i((Tf)_{(\gamma,0)})$.  Moreover, we let $k \in \{1,
  \dots, m\}$ and a $k$-zero $\kappa \in [0,\infty)^m$ be sufficiently
  small.  If $i \le k$, then $\t_{(\kappa,0)}(\partial_i f_{(\gamma,0)})
  = \partial_i(\t_{(\kappa,0)} f_{(\gamma,0)})$ by Remark
  \ref{translate_derivative}.  If $i > k$, then $\t_{(\kappa,0)}
  f_{(\gamma,0)}$ belongs to $\Q_k$ by Proposition \ref{Q-basics}(2),
  and it follows from Corollary \ref{good_derivatives} and Remark
  \ref{translate_derivative} that
  \begin{equation*}
    x_i^1 \cdot \frac {\partial}{\partial x_i} (\t_{(\kappa,0)}
    f_{(\gamma,0)}) = \partial_i (\t_{(\kappa,0)} f_{(\gamma,0)}) =
    \frac {x_i^1} {\kappa_i + x_i^1} \cdot \t_{(\kappa,0)} (\partial_i
    f_{(\gamma,0)}).
  \end{equation*}
  Therefore, $\t_{(\kappa,0)} (\partial_i f_{(\gamma,0)})$ belongs to
  $\A_m(V)$ for some suitable $V$ with $T(\t_{(\kappa,0)} (\partial_i
  f_{(\gamma,0)})) = (\kappa_i + X_i) \cdot (\partial/\partial
  X_i)(T(\t_{(\kappa,0)} f_{(\gamma,0)}))$, which proves part (1).

  Part (2) is more straightforward and follows from Corollary
  \ref{good_derivatives}, Proposition \ref{Q-basics} and Lemmas
  \ref{derivation_truncation} and \ref{Q-addition}.  For (3), note
  that for all $\gamma \in [0,\infty)^m$, we can take $g_{(\gamma,0)}
  := f_{(\gamma,0)}(x,y_1, \dots, y_{n-1},0)$.  Then for every
  sufficiently small $\kappa \in [0,\infty)^m$, we have
  $$\t_{(\kappa,0)} g_{(\gamma,0)} = (\t_{(\kappa,0)}
  f_{(\gamma,0)})(x,y_1, \dots, y_{n-1},0),$$ and part (3) follows.
\end{proof}

\subsection*{Composition}

Let $f \in \Q_m(U)$.  For the next lemma, we let $\q = (\q_1, \dots,
\q_n) \in \NN^n$ and put $\k := |\q|$.  We let $z = (z_1, \dots,
z_\k)$ range over $\frL^\k$ and $U' \subseteq \frL^{m+\k}$ be an
$m$-quadratic domain such that $(x, z_1 + \cdots + z_{\q_1}, \dots,
z_{\q_1 + \cdots + \q_{n-1}+1} + \cdots + z_\k) \in U$ for all $(x,z)
\in U'$.  In this situation, we define the holomorphic function $f_\q:
U' \into \CC$ by $$f_\q(x,z) := f(x, z_1 + \cdots + z_{\q_1}, \dots,
z_{\q_1 + \cdots + \q_{n-1}+1} + \cdots + z_\k).$$

\begin{lemma}
  \label{sum_composition}
  We have $f_\q \in \Q_m(U')$ and $T(f_\q) = (Tf)_\q$.
\end{lemma}

\begin{proof}
  We first show that $f_\q \in \A_m(U')$ and $T(f_\q) = (Tf)_\q$.
  Arguing by induction on $\k$ (simultaneously for all $m$) and
  permuting the last $n$ coordinates if necessary, it suffices to
  consider the case where $n = 1$ and $\k = \q_1 = 2$.  In this
  situation, by Proposition \ref{coefficients} and after shrinking $U$
  if necessary, we can write $f(x,y) = \sum_{p \in \NN} a_p(x) y^p$
  for all $(x,y) \in U$, and there are a quadratic domain $W \subseteq
  \frL$ and constants $A,B>0$ such that $a_p \in \A_m(W^m)$ and
  $\|a_p(x)\| \le AB^p$ for all $x \in W^m$ and each $p \in \NN$.
  Hence
  \begin{equation*}
    f_\q(x,z) = \sum_{p \in \NN} a_p(x)
    (z_1+z_2)^p = \sum_{p,q \in 
      \NN} b_{p,q}(x) z_1^p z_2^q
  \end{equation*}
  for all sufficiently small $(x,z) \in U'$, where $b_{p,q} :=
  \begin{pmatrix} p+q \\ p \end{pmatrix} a_{p+q}$ for all $p,q \in
  \NN$.  Since $\|b_{p,q}\| \le A (2B)^{p+q}$, it follows from
  Propositions \ref{coefficients} and \ref{O-membership} that $f_\q
  \in \A_m(U')$, as required.

  Next, for every $\gamma \in [0,\infty)^m$, we have
  $T(f_{(\gamma,0)})_\q) = (T f_\q)_{(\gamma,0)}$, so we can take
  $(f_\q)_{(\gamma,0)} := (f_{(\gamma,0)})_\q$.  Moreover, for every
  sufficiently small $\kappa \in [0,\infty)^m$, the previous paragraph
  and Proposition \ref{O-membership} now also show that
  $\t_{(\kappa,0)} (f_\q)_{(\gamma,0)}$ belongs to $\A_m(V')$ for some
  appropriate $V'$.
\end{proof}

For the next proposition, we let $g = (g_1, \dots, g_n) \in \Q_m(V)^n$
be such that $g(0) = 0$ and $(x, g(x,y)) \in U$ for all $(x,y) \in V$,
and we define the holomorphic function $h:V \into \CC$ by $ h(x,y) :=
f(x,g(x,y)).  $

\begin{prop}
  \label{Q-composition}
  The function $h$ belongs to $\Q_m(V)$ and $Th(X,Y) =
  Tf(X,Tg(X,Y))$.
\end{prop}

\begin{proof}
  First, let $\kappa \in [0,\infty)^m$ be such that $(\kappa,0) \in
  \cl_0(V)$.  Then $\t_{(\kappa,0)} h =
  (\t_{(\kappa,0)}f)(x,\t_{(\kappa,0)} g)$, so Corollary \ref{taylor}
  (with $\lambda$ there equal to $\t_{(\kappa,0)}g(0,0)$) and
  Proposition \ref{A-composition} show that $\t_{(\kappa,0)}h \in
  \A_m(W)$ for some appropriate $W$.

  Second, let $\gamma \in [0,\infty)^m$; we need to find an
  $m$-quadratic domain $V' \subseteq V$ and an $h' \in \A_m(V')$ such
  that $T(h') = (Th)_{(\gamma,0)}$.  By Proposition
  \ref{comp_representation}, there are $p \in \NN$, a tuple $\q \in
  \NN^n$ and, with $\k:= |\q|$, elementary sets $E_1, \dots, E_p
  \subseteq [0,\infty)^m \times \NN^\k$ and $B_{i,j} \subseteq
  [0,\infty)^m \times \NN^n$ for each pair $(i,j)$ satisfying $i \in
  \{1, \dots, n\}$ and $j \in \{1, \dots, \q_i\}$, such that
  \begin{equation*}
    Tf(X,Tg)_{(\gamma,0)} = \sum_{q=1}^p \frac {(X,(Tg)_B)^{\inf
        E_q}}{X^\gamma} \cdot ((Tf)_\q)_{E_q}(X,(Tg)_B) 
  \end{equation*}
  with $(Tg)_B := ((Tg_1)_{B_{1,1}}, \dots, (Tg_n)_{B_{n,\q_n}})$ and
  each $(X,(Tg)_B)^{\inf E_q}$ divisible by $X^\gamma$.  After
  shrinking $V$ if necessary and writing $g_B := ((g_1)_{B_{1,1}},
  \dots, (g_n)_{B_{n,\q_n}})$, we get from Lemma
  \ref{sum_composition}, Proposition \ref{Q-basics} and the above
  that, for each $q =1 , \dots, p$, the function $$h_q := x^{-\gamma}
  \cdot (x,g_B)^{\inf E_q} \cdot (f_\q)_{E_q}(x,g_B)$$ belongs to
  $\A_m(V)$ and satisfies $$T h_\q = X^{-\gamma} \cdot
  (X,(Tg)_B)^{\inf E_q} \cdot ((Tf)_\q)_{E_q}(X,(Tg)_B).$$ Hence by
  Lemma \ref{Q-addition}, we can take $h' := h_1 + \cdots + h_p$.

  Finally, it follows from the last paragraph and the first
  observation above that $h \in \Q_m$.
\end{proof}

Here are some immediate applications of Proposition \ref{Q-composition}:

\begin{prop}
  \label{Q-algebra}
  The set $\Q_m$ is a $\CC$-algebra, and the map $f \mapsto Tf: \Q_m
  \longrightarrow \Ps{C}{X^*,Y}$ is an injective $\CC$-algebra
  homomorphism such that $f(0) = (Tf)(0)$ for all $f \in \Q_m$.
\end{prop}

\begin{proof}
  Let $f,g \in \Q_m$; we need to show that $fg \in \Q_m$.  Put $f_1 :=
  f - f(0)$ and $g_1 := g - g(0)$; then $f_1, g_1 \in \Q_m$ by Lemma
  \ref{Q-addition}, and $fg = P(f_1,g_1)$ with $P(Y_1,Y_2) := (f(0) +
  Y_1) (g(0) + Y_2)$.  Hence $fg \in \Q_m$ by Proposition
  \ref{Q-composition}.
\end{proof}

\begin{prop}
  \label{Q-units}
  Let $f \in \Q_m$.  Then 
  \begin{enumerate}
  \item $f$ is a unit in $\Q_m$ if and only if $f(0) \ne 0$;
  \item if $f(0) = 0$, then there are $(\gamma,\delta) \in
    (0,\infty)^m \times (\NN \setminus \{0\})^n$ and $f_1, \dots,
    f_{m+n} \in \Q_m$ such that $f = X_1^{\gamma_1} f_1 + \cdots +
    Y_n^{\delta_n} f_{m+n}$.
  \end{enumerate}
\end{prop}

\begin{proof}
  (1) Assume first that $f$ is a unit in $\Q_m$, and let $g \in Q_m$
  be such that $f\cdot g = 1$.  Then $Tf \cdot Th = 1$ and hence $f(0)
  = Tf(0) \ne 0$.  Conversely, assume that $f(0) \ne 0$; we may assume
  that $f(0) = 1$, and we put $f_1:= 1 - f \in \Q_m$.  Let $\a>0$ and
  $U_\a \subseteq \frL^{m+n}$ be an $m$-quadratic domain such that
  $f_1(x,y) = o(\|(x,y)\|^{\a})$ as $\|(x,y)\| \to 0$ in $U_{\a}$.
  Thus, there is an $m$-quadratic domain $U \subseteq U_{\a}$ such
  that $\|f_1(x,y)\| \leq \frac{1}{2}$ for all $(x,y) \in U$.  Let
  $\phi: B(0,1) \into \CC$ be the holomorphic function defined by
  $\phi(z) := \frac{1}{1-z}$, and define $g:U \into \CC$ by $g(x,y) :=
  \phi(f_1(x,y))$.  Then $f \cdot g = 1$, and $g \in \Q_m(U)$ by
  Proposition \ref{Q-composition}.

  (2) follows from Lemma 4.8 of \cite{vdd-spe:genpower} and
  Proposition \ref{Q-basics}.
\end{proof}

Finally, Proposition \ref{Q-composition} allows us to make sense of
certain substitutions in the $x$-variables:

\begin{df}
  \label{x-substitution}
  Let $W \subseteq \frL$ be a quadratic domain and $R>0$, and let $f
  \in \Q_m(W^m \times B_\frL(R))$.  Let also $V \subseteq \frL^{m+n}$ be
  $m$-quadratic and $g = (g_1, \dots, g_m) \in \Q_m(V)^m$ be such that
  $\lambda := g(0,0) \in W^m \cap (0,\infty)^m$.  Then $g(x,y) =
  \lambda + h(x,y)$ with $h \in \Q_m(V)^m$ satisfying $h(0) = 0$, and
  we define $f(g(x,y),y) := (\t_{(\lambda,0)} f) (h(x,y),y)$.
\end{df}

\begin{cor}
  \label{bad_substitution_1}
  The function $f(g(x,y),y)$ in Definition \ref{x-substitution}
  belongs to $\Q_m$. \qed
\end{cor}

Some of the substitutions not covered by the previous corollary are
the blow-up substitutions:

\begin{prop}
  \label{blowup_homom}
  Let $\rho,\lambda > 0$ and $i,j \in \{1, \dots, m\}$ be distinct.
  \begin{enumerate}
  \item The function $f \circ \s^\rho_{i,j}$ belongs to $\Q_m$ for
    every $f \in \Q_m$, and the map $\s^\rho_{ij}:\Q_m \into \Q_m$
    defined by $\s^\rho_{ij}(f) := f \circ \s^\rho_{ij}$ is a
    $\CC$-algebra homomorphism such that $T \circ \s^\rho_{ij} =
    \frB^{\rho,0}_{ij} \circ T$.
  \item The function $f \circ \r^{\rho,\lambda}$ belong to $\Q_{m-1}$
    for every $f \in \Q_m$, and the map $\r^{\rho,\lambda}:\Q_m \into
    \Q_{m-1}$ defined by $\r^{\rho,\lambda}(f) := f \circ
    \r^{\rho,\lambda}$ is a $\CC$-algebra homomorphism such that $T
    \circ \r^{\rho,\lambda} = \frB^{\rho,\lambda}_{m,m-1} \circ T$.
  \end{enumerate}
\end{prop}

Whenever convenient, we shall write $\s^\rho_{ij} f$ and
$\r^{\rho,\lambda} f$ in place of $\s^\rho_{ij}(f)$ and
$\r^{\rho,\lambda}(f)$.

\begin{proof}
  The proofs for parts (1) and (2) are similar; we prove (1) here and
  leave (2) to the reader.  We may assume that $i=m$ and $j=m-1$, and
  we write $\s$ and $\B$ in place of $\s^\rho_{m,m-1}$ and
  $\B^\rho_{m,m-1}$.  Let $f \in \Q_m$; if suffices to prove that $f
  \circ \s \in \Q_m$.  

  To do so, we let $W \subseteq \frL$ be quadratic and $1 > R>0$ be
  such that $f \in \A_m(W^m \times B_\frL(R)^n)$, and we let $W'$ and $V$
  be as in the proof of Proposition \ref{singular_prop}.  We also let
  $\kappa \in [0,\infty)^m$ be nonzero such that $(\kappa,0) \in
  \cl_0(V)$, and let $W_\kappa \subseteq W'$ be quadratic such that
  $\t_{(\kappa,0)}(x,y) \in V$ for all $(x,y) \in V_\kappa :=
  (W_\kappa)^m \times B_\frL(R)^n$.  By Propositions
  \ref{blowup_truncation} and \ref{Q-basics}, it remains to prove that
  \begin{itemize}
  \item[$(\ast)$] $\t_{(\kappa,0)}(f \circ \s)$ belongs to
    $\A_m(V_\kappa)$.  
  \end{itemize}
  Writing $\kappa'\:= (\kappa_1, \dots, \kappa_{m-2})$ and $\kappa'':=
  (\kappa_{m-1}, \kappa_m)$, we see that $\t_{(\kappa,0)}(f \circ \s)
  = \t_{(0,\kappa'',0)}( \t_{(\kappa',0,0)}(f \circ \s))$; since
  $\t_{(\kappa',0,0)}(f \circ \s) = (\t_{(\kappa',0,0)} f) \circ \s$,
  we may even assume that $\kappa_1 = \cdots = \kappa_{m-2} = 0$.  We
  now distinguish three cases: \medskip

  \noindent\textbf{Case 1:} both $\kappa_{m-1}$ and $\kappa_m$ are
  nonzero.  Then
  \begin{equation*}
    \t_{(\kappa,0)}(f \circ \s)(x,y) =
    (\t_{(0,\kappa_{m-1},\kappa_{m-1}^\rho \kappa_m,0)}
    f)(x',g(x_{m-1}, x_m),y),
  \end{equation*}
  where $x' := (x_1, \dots, x_{m-1})$ and $g$ is an analytic function
  satisfying $g(0) = 0$.  Since $\t_{(0,\kappa_{m-1},\kappa_{m-1}^\rho
    \kappa_m,0)} f$ belongs to $\Q_{m-2}$ by Proposition
  \ref{Q-basics}, $(\ast)$ follows from Proposition
  \ref{Q-composition} in this case.  \medskip

  \noindent\textbf{Case 2:} $\kappa_{m-1} = 0$ and $\kappa_m > 0$.
  Then $\t_{(\kappa,0)} (f \circ \s) = f \circ \r^{\rho,\kappa_m}$,
  so $(\ast)$ follows from Proposition \ref{regular_prop} in this
  case. 
  \medskip

  \noindent\textbf{Case 3:} $\kappa_{m-1} \ne 0$ and $\kappa_m = 0$.
  We define $\phi(z_1, \dots, z_{m+1},y) := \t_{(\kappa,0)} f(z_1,
  \dots, z_{m-1},z_{m+1},y)$.  Then $\phi \in \Q_{m+1}$, and there is
  an analytic one-variable function $g$ with $g(0) = 0$ such that
  \begin{equation*}
    \t_{(\kappa,0)} (f \circ \s)(x,y) = \left(\phi \circ
    \r^{1,\kappa_{m-1}^\rho}\right) (x, g(x_{m-1}),y).  
  \end{equation*}
  Thus, $(\ast)$ follows from Propositions \ref{regular_prop} and
  \ref{A-composition} in this case.
\end{proof}

As a consequence of Proposition \ref{blowup_homom}, we extend
Corollary \ref{bad_substitution_1} to certain functions with zero
constant coefficient:

\begin{df}
  \label{x-substitution_2}
  Let $m \ge 1$, let $W \subseteq \frL$ be a quadratic domain and
  $R>0$, and let $f \in \Q_m(W^m \times B_\frL(R)^n)$.  Let also $V
  \subseteq \frL$ be a quadratic domain and $g \in Q_1(V)$ be such
  that $g(t) \in W \cap (0,\infty)$ for all $t \in V \cap (0,\infty)$.
  Then $g(t) = t^\rho(\lambda + h(t))$ for some $\rho, \lambda > 0$
  and some $h \in \Q_1(V)$ with $h(0) = 0$.  We write $x' := (x_1,
  \dots, x_{m-1})$ and let $\tilde{f} \in \Q_{m+1}(W^{m+1} \times
  B_\frL(R)^n)$ be the function defined by $\tilde{f}(x',u,v,y) :=
  f(x',v,y)$.  Then $\r^{\rho,\lambda} \tilde{f} \in \Q_m^{m+n+1}$,
  and we define
  $$f(x',g(t),y) := \left(\r^{\rho,\lambda} \tilde{f}\right)
  (x',t,h(t),y).$$
\end{df}

\begin{cor}
  \label{bad_substitution_2}
  The function $f(x',g(t),y)$ in Definition \ref{x-substitution_2}
  belongs to $\Q_m^{m+n}$. \qed
\end{cor}

\section{Weierstrass Preparation}  \label{weierstrass}

We continue to work in the setting of the previous section.  In this
section, we establish a Weierstrass Preparation Theorem for the
classes $\Q_m$.  We follow Brieskorn and Kn\"orrer's exposition in
Section 8.2 of \cite{bri-kno:plane}; to do so, we need to first
establish an implicit function theorem and a theorem on symmetric
functions.  We thank Lou van den Dries for his helpfull suggestions on
this section, especially the proof of Corollary \ref{multi_implicit}
below.

We start with a single implicit variable and write $y' :=
(y_1, \dots, y_{n-1})$ and $Y' := (Y_1, \dots, Y_{n-1})$.

\begin{prop}
  \label{implicit}
  Let $f \in \Q_m^{m+n}$, and assume that $f(0) = 0$ and $\partial
  f/\partial y_n(0) \ne 0$.  Then there is an $h \in \Q_m^{m+n-1}$ such
  that $h(0) = 0$ and $f(x,y',h) = 0$.
\end{prop}

\begin{proof}
  We let $U = W^m \times B_\frL(R)^n$ for some quadratic $W \subseteq \frL$
  and $R>0$ be such that $f \in \Q_m(U)$, and we put $V := W^m \times
  B_\frL(R)^{n-1}$.  During the proof below, we may have to shrink $W$ and
  $R$ (and all related quantities introduced below) on various
  occasions; we will not explicitely mention this.  By Proposition
  \ref{Q-derivatives}, the function $(\partial f/\partial
  y_n)(x,y',0)$ belongs to $\Q_m(V)$.  Hence by hypothesis, there is a
  constant $c>0$ such that $|(\partial f/\partial y_n)(x,y',0)| \ge c$
  for all $(x,y') \in V$.  On the other hand, by Proposition
  \ref{coefficients}, we can write $f(x,y) = \sum_{p \in \NN} a_{p}(x,y')
  y_n^p$ with each $a_{p} \in \O_m(V)$.

  Define $g:U \into \CC$ by $g(x,y):= f(x,y) - a_0(x,y')$; note that
  $g \in \O_m(U)$.  By the above and the usual arguments for the
  inverse function theorem, there is a $\rho>0$ such that for every
  $(x,y') \in V$, we have $\|a_0(x,y')\| \le \rho/2$ and the function
  $g_{x,y'}:B_\frL(R) \into \CC$ defined by $g_{x,y'}(y_n) := g(x,y)$ is
  injective and satisfies $g_{x,y'}(0) = 0$ and $B(0,\rho) \subseteq
  g_{x,y'}(B_\frL(R))$, and such that its compositional inverse
  $g_{x,y'}^{-1}:B_\frL(\rho) \into \CC$ is given by a convergent power
  series
  \begin{equation*}
    g^{-1}_{x,y'}(z) = \sum_{p \in \NN} b_p(x,y') y_n^p.
  \end{equation*}

  We claim that the function $H:V \times B_\frL(\rho/2) \into \CC$ defined
  by $H(x,y) := g_{x,y'}^{-1}(y_n)$ belongs to $\Q_m(V \times
  B_\frL(\rho/2))$.  The proposition follows from this claim by defining
  $h:V \into \CC$ as $h(x,y') := H(x,y',y_n-a_0(x,y'))\rest{y_n=0}$.

  To see the claim, we note first from the Lagrange inversion formula
  (see for instance Whittaker and Watson
  \cite[p. 133]{whi-wat:analysis}) that for all $p \in \NN$,
  \begin{equation*}
    b_p(x,y') = \frac1{p!} \left[ \frac{\partial^{p-1}}{\partial y_n^{p-1}}
      \left( \frac {y_n}{g(x,y)} \right)^p \right]_{y_n=0}.
  \end{equation*}
  Since $a_1(0) \ne 0$ by hypothesis, it follows from Propositions
  \ref{Q-units}, \ref{Q-algebra} and \ref{Q-derivatives} and Remark
  \ref{quasirmks}(2) that $b_p \in \Q_m(V)$ for each $p$.

  Second, we let $\gamma \in [0,\infty)^m$, and we claim that
  $H_{(\gamma,0)} \in \A_m(V \times B_\frL(\rho/2))$.  It suffices, by
  Propositions \ref{O-membership} and \ref{holo_division}, to find
  constants $A,B > 0$ such that $\|(b_p)_{(\gamma,0)}(x,y')\| \le A
  B^p$ for all $p \in \NN$ and $(x,y') \in V$.  To do so, we shall
  assume that $U \subseteq B_\frL(1)^{m+n}$, and we put $\g(x,y) :=
  y_n/g(x,y) \in \Q_m(U)$.  We let $\J$ be the set of all ordered
  pairs $(I,\alpha)$ such that $I \subseteq \{1, \dots, m\}$ and
  $\alpha \in B_I = B_{\gamma,I}(T\mathfrak{g})$, where the latter is
  defined as in Lemma \ref{lemma:representation}.  By Corollary
  \ref{cor:Q-basics}, there is a constant $C>0$ and for each
  $(I,\alpha) \in \J$ a function $\g_{I,\alpha} \in \Q_m(U)$,
  depending only on the variables $x_{\bar I}$ and $y$, such that
  \begin{equation*}
    \g(x,y) = \sum_{(I,\alpha) \in \J} x_{\bar I}^{\gamma_{\bar I}}
    x_I^\alpha \g_{I,\alpha}(x,y)
  \end{equation*}
  and $\|\g_{I,\alpha}(x,y)\| \le C$ for all $(x,y) \in U$.  We fix $p
  \in \NN$ and write $\g^p(x,y) := (\g(x,y))^p$ for all $(x,y) \in U$.
  Since
  \begin{multline*}
    \g^p(x,y) = \\ \sum_{((I_1,\alpha_1), \dots, (I_p, \alpha_p))
      \in \J^p} x_{I_1}^{\alpha_1} \cdots x_{I_p}^{\alpha_p} \cdot
    x_{\bar{I_1}}^{\gamma_{\bar{I_1}}} \cdots
      x_{\bar{I_p}}^{\gamma_{\bar{I_p}}} \cdot
        \g_{I_1,\alpha_1}(x,y) \cdots
        \g_{I_p,\alpha_p}(x,y),
  \end{multline*}
  and since each $\mathfrak{g}_{I,\alpha}$ only depends on the
  variables $x_{\bar I}$ and $y$, we get that
  \begin{multline*}
    (\g^p)_{(\gamma,0)}(x,y) = \\ \sum_{((I_1,\alpha_1), \dots, (I_p,
      \alpha_p)) \in \J_p} \frac {x_{I_1}^{\alpha_1} \cdots
      x_{I_p}^{\alpha_p} \cdot x_{\bar{I_1}}^{\gamma_{\bar{I_1}}}
      \cdots x_{\bar{I_p}}^{\gamma_{\bar{I_p}}}} {x^\gamma} \cdot
    \g_{I_1,\alpha_1}(x,y) \cdots
    \g_{I_p,\alpha_p}(x,y),
  \end{multline*}
  where $\J_p$ is the set of all $((I_1, \alpha_1), \dots,
  (I_p,\alpha_p)) \in \J^p$ such that the monomial $x_{I_1}^{\alpha_1}
  \cdots x_{I_p}^{\alpha_p} \cdot x_{\bar{I_1}}^{\gamma_{\bar{I_1}}}
  \cdots x_{\bar{I_p}}^{\gamma_{\bar{I_p}}}$ is divisible by
  $x^\gamma$.  As $U \subseteq B_\frL(1)^{m+n}$, it follows for all $x
  \in U$ that
  \begin{align*}
    \left\|(\g^p)_{(\gamma,0)}(x,y)\right\| &\le
    \sum_{((I_1,\alpha_1), \dots, (I_p, \alpha_p)) \in \J^p}
    \|\g_{I_1,\alpha_1}(x,y)\| \cdots
    \|\g_{I_p,\alpha_p}(x,y)\| \\ &= \left(\sum_{(I,\alpha) \in \J}
      \|\g_{I,\alpha}(x,y)\|\right)^p \\ &\le \left(|\J|
      C\right)^p. 
  \end{align*}
  It follows from Lemma \ref{derivation_truncation} and the Cauchy
  estimates that
  \begin{equation*}
    \|(b_p)_{(\gamma,0)}(x,y')\| = \frac 1{p!} \left\| \frac {\partial^{p-1}
        (\g^p)_{(\gamma,0)}} {\partial y_n^{p-1}} (x,y',0) \right\| \le \frac
    {\left(|\J|C\right)^p} {\rho^{p-1}} 
  \end{equation*}
  for all $(x,y') \in V$; so we can take $A := \rho$ and $B :=
  |\J|C/\rho$.

  Finally, let $\kappa \in [0,\infty)^m$ be such that $(\kappa,0) \in
  \cl_0(V \times B_\frL(\rho/2))$.  Then $\t_{(\kappa,0)} (b_p)_{(\gamma,0)}
  \in \A_m(V')$ for some appropriate $V'$ independent of $p$ by Remark
  \ref{O-rmk}(2).  Since $\|\t_{(\kappa,0)} (b_p)_{(\gamma,0)}(x,y')\|
  \le AB^p$ for all $(x,y') \in V'$ by the above, it follows that
  $\t_{(\kappa,0)} H_{(\gamma,0)}$ belongs to $\A_m(V' \times
  B_\frL(\rho/2))$, and the proposition is proved.
\end{proof}

The case of several implicit variables can be reduced to that of one
implicit variable: below, we let $l \in \{1, \dots, n\}$, and we write
$y' := (y_1, \dots, y_{n-l})$, $z = (z_1, \dots, z_l):= (y_{n-l+1},
\dots, y_n)$, $Y' := (Y_1, \dots, Y_{n-l})$ and $Z = (Z_1, \dots, Z_l)
:= (Y_{n-l+1}, \dots, Y_n)$.

\begin{cor}[Implicit Function Theorem]
  \label{multi_implicit}
  Let $f \in (\Q_m)^l$  such that $f(0) = 0$ and $\partial
  f/\partial z(0) \ne 0$.  Then there is an $h \in (\Q_m^{m+n-l})^l$
  such that $h(0) = 0$ and $f(x,y',h) = 0$.
\end{cor}

\begin{proof}
  By induction on $l$; the case $l=1$ corresponds to Proposition
  \ref{implicit}, so we assume that $n \ge l>1$ and the corollary
  holds for lower values of $l$.  After permuting the component
  functions of $f$, we may assume that $\partial f_l/\partial y_n(0)
  \ne 0$.  Hence, writing $z' := (y_{n-l+1}, \dots, y_{n-1})$, we
  obtain from Proposition \ref{implicit} a function $w \in
  \Q_m^{m+n-1}$ such that $f_l(x,y', z', w) = 0$.  Moreover, there are
  constants $c_1, \dots, c_{l-1} \in \CC$ such that for each $i=1,
  \dots, l-1$, defining $f'_i := f_i - c_i f_l$ gives $(\partial
  f'_i/\partial y_n)(0) = 0$.  By the hypothesis of the corollary, the
  map $g \in (\Q_m^{m+n-1})^{l-1}$ defined by
  \begin{equation*}
    g_i(x,y',z') := f'_i(x,y',z',w(x,y',z')) \quad\text{for } i=1,
    \dots, l-1 
  \end{equation*}
  satisfies $g(0) = 0$ and $(\partial g / \partial z')(0) \ne 0$.
  Hence by the inductive hypothesis, there is an $h' \in
  (\Q_m^{m+n-l})^{l-1}$ such that $g(x,y',h') = 0$.  The corollary
  follows with $h \in (\Q_m^{m+n-l})^l$ defined by $h_i := h'_i$ if $i
  = 1, \dots, l-1$ and $h_l(x,y') := w(x,y',h'(x,y'))$.
\end{proof}

For the next proposition, we let $\sigma = (\sigma_1, \dots,
\sigma_l)$ be the elementary symmetric functions in the variables $z$.
Recall that $f \in \Q_m$ is \textbf{symmetric in the variables $z$} if
$f(x,y',z) = f(x,y', \lambda(z))$ for every permutation $\lambda$ of
$\{1, \dots, l\}$.

\begin{prop}[Symmetric Function Theorem]
  \label{symmetric}
  Let $f \in \Q_m$ be symmetric in the variables $z$.  Then there is a
  $g \in \Q_m$ such that $f(x,y',z) = g(x,y',\sigma)$.
\end{prop}

\begin{proof}
  First, let $\gamma \in [0,\infty)^m$ and assume there is a $G \in
  \Ps{C}{X^*,Y',Z}$ such that $Tf(X,Y',Z) = G(X,Y',\sigma_1(Z), \dots,
  \sigma_l(Z))$.  Then $(Tf)_{(\gamma,0,0)}$ is symmetric in $Z$
  and $$(Tf)_{(\gamma,0,0)}(X,Y',Z) =
  G_{(\gamma,0,0)}(X,Y',\sigma_1(Z), \dots, \sigma_l(Z)).$$ Moreover,
  if $g \in \A_m(V)$ is such that $f(x,y',z) = g(x,y',\sigma)$, and if
  $\kappa \in [0,\infty)^m$ is sufficiently small, then
  $\t_{(\kappa,0,0)} f$ is symmetric in $z$ and $\t_{(\kappa,0,0)}
  f(x,y',z) = \t_{(\kappa,0,0)} g(x,y',\sigma)$.

  Therefore, we assume that $f \in \A_m(U)$ for some $m$-quadratic $U
  \subseteq \frL^{m+n}$ and we need to find, after shrinking $U$ if
  necessary, a $g \in \A_m(U)$ such that $f(x,y',z) = g(x,y',\sigma_1,
  \dots, \sigma_l)$ for all $(x,y',z) \in U$.  Without loss of
  generality, we also assume that $U = W^m \times B_\frL(R)^n$ for some
  quadratic domain $W \subseteq \frL$ and some $R>0$, and we put $U'
  := W^m \times B_\frL(R)^{n-l}$.

  By Proposition \ref{coefficients}, there are $a_q \in
  \A_m^{m+n-l}(U')$, for $q \in \NN^l$, and constants $B,C>0$ such
  that $\|a_q(x,y')\| \le B C^{|q|}$ for each $q$ and $f(x,y',z) =
  \sum_{q \in \NN^l} a_q(x,y') z^q$.  Let also $\sim$ be the
  equivalence relation on $\NN^l$ defined by $p \sim q$ if and only if
  there is a permutation $\lambda$ of $\{1, \dots, l\}$ such that $p =
  (q_{\sigma(1)}, \dots, q_{\sigma(l)})$, and let $E_1, E_2, \dots$ be
  an enumeration of all equivalence classes of $\sim$.  Since $f$ is
  symmetric in $z$, we get for all $j \in \NN$ that $a_p = a_q$ for
  all $p,q \in E_j$.  Thus, for each $j \in \NN$,  we define $b_j :=
  a_p$ for some $p \in E_j$; then
  \begin{equation*}
    \sum_{p \in E_j} a_p(x,y') z^p = b_j(x,y') \cdot \sum_{p \in E_j}
    z^p \quad\text{for all } j \in \NN.
  \end{equation*}
  Let $j \in \NN$, and note that the sum $\sum_{p \in E_j} z^p$ is a
  symmetric polynomial in $z$ that is homogeneous of degree $d_j :=
  |p|$ for any $p \in E_j$.  By the main theorem on symmetric
  polynomials (see for instance Van der Waerden \cite{vdw:algebra}),
  there is a unique polynomial $S_j \in \CC[Z]$ of weighted degree
  $d_j$ such that $\sum_{p \in E_j} z^p = S_j(\sigma_1(z), \dots,
  \sigma_l(z))$.  (Here ``of weighted degree $d_j$'' means that any
  term $c z^p$ occurring in $S_j$ satisfies $p_1 + 2p_2 + \cdots +
  lp_l = d_j$.)

  We now define $S \in \Ps{C}{Z}$ by $S(Z) := \sum_{j \in \NN}
  S_j(Z)$; we claim that $S$ is convergent.  To see this, note that
  the product $\Pi_{i=1}^l (1-Z_i)$ is a symmetric polynomial; hence,
  there is a polynomial $P \in \CC[Z]$ such that $P(\sigma_1(Z),
  \dots, \sigma_l(Z)) = \Pi_{i=1}^l (1-Z_i)$.  Since $P(0) \ne 0$, we
  have that $1/P \in \Ps{C}{Z}$ converges; but $(1/P)(\sigma_1(Z),
  \dots, \sigma_l(Z)) = \sum_{q \in \NN^l} z^q$, and the claim
  follows. 

  We now claim that the sum $g(x,y',z) := \sum_{j \in \NN} b_j(x,y')
  S_j(z)$ defines a function $g \in \A_m(U)$.  Assuming the claim
  holds, we necessarily have $Tg(X,Y',Z) = G(X,Y',Z) := \sum_{j \in \NN}
  Tb_j(X,Y') S_j(Z)$, so $g$ has the required properties by the
  injectivity of the map $T$ and because $Tf(X,Y',Z) = G(X,Y',
  \sigma_1(Z), \dots, \sigma_l(Z))$.

  To see the claim, we first need to rewrite the sum: for each $q \in
  \NN^l$, we put $D_q := \{j \in \NN:\ q \in \supp S_j\}$ and $c_q :=
  \sum_{j \in D_q} b_j$, so that $g(x,y',z) = \sum_{q \in \NN^l}
  c_q(x,y') z^q$ for all $(x,y',z) \in U$.  Note that for all $q \in
  \NN^l$, $j \in D_q$ and $p \in E_j$, we have $|p| = q_1 + 2q_2 +
  \cdots + lq_l \le l|q|$.  Hence, for each $q \in \NN^l$, there is a
  set $C_q \subseteq \set{r \in \NN^l:\ |r| \le l|q|}$ such that $c_q
  = \sum_{r \in C_q} a_r$.  Since $|C_q| \le (l|q|)^l$, it follows
  that there are constants $\a,\b > 0$ such that $\|c_q(x,y')\| \le
  \a\b^{|q|}$ for every $q \in \NN^l$ and all $(x,y') \in U'$.  Since
  $c_q \in \A_m^{m+n-l}(U')$ for every $q \in \NN^l$, the claim
  follows from Proposition \ref{O-membership}.
\end{proof}

For the remainder of this section, we write again $Y' = (Y_1, \dots,
Y_{n-1})$.  We recall from Definition 4.16 of \cite{vdd-spe:genpower}
that $F \in \Ps{C}{X^*,Y}$ is regular in $Y_n$ of order $d \in \NN$ if
$F(0,0,Y_n) = c Y_n^d + $ terms of higher order in $Y_n$.  

\begin{prop}[Weierstrass Preparation]
  \label{wp_theorem}
  Assume that $n>0$, and let $f \in \Q_m$ be such that $Tf$ is regular
  in $Y_n$ of order $d \in \NN$.  
  \begin{enumerate}
  \item For every $g \in \Q_m$, there are a unique $q \in \Q_m$ and a
    unique $r \in \Q_m^{m+n-1}[Y_n]$ such that $g = qf + r$ and
    $\deg_{Y_n}(r) < d$.
  \item There are a unique unit $u \in \Q_m$ and a unique $w \in
    \Q_m^{m+n-1}[Y_n]$ such that $f = uw$ and $w$ is monic of degree
  $d$ in $Y_n$.
  \end{enumerate}
\end{prop}

\begin{proof}
  The proof of Theorems 1 and 2 on p. 338 of \cite{bri-kno:plane} goes
  through almost literally, using the properties established for the
  classes $\Q_m$ in the previous sections as well as the Implicit
  Function Theorem and the Symmetric Function Theorem above, except
  for the following trivial changes: the variable $t$ and $z_1, \dots,
  z_n$ there correspond to $y_n$ and $x_1, \dots, x_m, y_1, \dots,
  y_{n-1}$ here, and the roles of $f$ and $g$ are exchanged.  (Note
  that the uniqueness also follows directly from Proposition 4.17 in
  \cite{vdd-spe:genpower} and the injectivity of the map $T:\Q_m \into
  \Ps{C}{X^*,Y}$.)
\end{proof}

\section{$\Q$-semianalytic sets and model completeness}
  \label{o-minimal}

In this section we prove model completeness and o-minimality of $\Rq$.
We also show that $\Rq$ admits analytic cell decomposition.  (An
o-minimal expansion $\Rtilde$ of the ordered field of real numbers is
said to admit \textbf{analytic cell decomposition} if for any $A_{1},
\dots, A_{k} \subseteq \RR^{m}$ definable in $\Rtilde$, there is a
decomposition of $\RR^{m}$ into analytic cells definable in $\Rtilde$
and compatible with each $A_{i}$.)

We let $m,n \in \NN$ and $\rho \in (0,\infty)^{m+n}$, and we put
\begin{equation*}
  I_{m,n;\rho} := [0,\rho_1] \times \cdots \times [0,\rho_m] \times
  [-\rho_{m+1}, \rho_{m+1}] \times \cdots \times [-\rho_{m+n},
  \rho_{m+n}],
\end{equation*}
a subset of $\RR^{m+n}$.  We also write $I_{m,n;\epsilon}$ instead of
$I_{m,n;(\epsilon, \dots, \epsilon)}$, for $\epsilon > 0$, and we put
$I_{m,n;\infty} := [0,\infty)^m \times \RR^n$.  Abusing notation, we
identify $I_{m,n;\rho}$ with the set
\begin{multline*}
  \{(x,y) \in [0,\infty)_{\frL_0}^m \times \RR^n:\ 0 \le \|x_i\| \le
  \rho_i \text{ and } \\ -\rho_{m+j} \le y_j \le \rho_{m+j} \
  \text{for } i=1, \dots, m \text{ and } j=1, \dots, n\}.
\end{multline*}
Given an $m$-quadratic $U \subseteq \frL^{m+n}$ such that
$I_{m,n;\rho} \subseteq \ir(\cl(\pi^{m+n}_m(U))$, and given an $f \in
\Q_m(U)$, we write $f\rest{I_{m,n;\rho}}$ for the function
$f^\sharp\rest{I_{m,n;\rho}}$.

\begin{df}
  \label{real_Q}
  We let $\Q_{m,n;\rho}$ be the set of all functions $f:I_{m,n;\rho}
  \into \RR$ for which there exist an $m$-quadratic domain $U
  \subseteq \frL^{m+n}$ and a $g \in \Q_m(U)$ such that $I_{m,n;\rho}
  \subseteq \ir(\cl(\pi^{m+n}_m(U))$ and $f = g \rest{I_{m,n;\rho}}$.
\end{df}

\begin{nrmk}
  \label{restricted}
  For every $f \in \Q_{m,n;\rho}$, there are $\rho' > \rho$ and $g \in
  \Q_{m,n,\rho'}$ such that $f = g\rest{I_{m,n;\rho}}$.
\end{nrmk}

\begin{prop}
  \label{specialize}
  Let $U \subseteq \frL^{m+n}$ be $m$-quadratic such that $I_{m,n;\rho}
  \subseteq \ir(\cl(\pi^{m+n}_m(U))$, and let $f \in \Q_m(U)$.  Then
  $f\rest{I_{m,n;\rho}} \in \Q_{m,n;\rho}$ if and only if $Tf \in
    \Ps{R}{X^*,Y}$. 
\end{prop}

\begin{proof}
  The necessity is clear from Definition \ref{asymptotic}, so we
  assume that $Tf \in \Ps{R}{X^*,Y}$.  Define $\bar{(r,\varphi)} :=
  (r,-\varphi)$ for $(r,\varphi) \in \frL$ and $\bar{(x,y)} :=
  (\bar{x_1}, \dots, \bar{y_n})$ for $(x,y) \in \frL^{m+n}$.  Then the
  function $\bar{f}:U \into \CC$ defined by $\bar{f}(x) :=
  \bar{f(\bar{x})}$ belongs to $\A(U)$ and satisfies
  $T\left(\bar{f}\right) = Tf$.  Hence by Proposition
  \ref{asymptotic_to_0}, we have $\bar{f} = f$, which proves the
  proposition.
\end{proof}

Correspondingly, we put $\Pc{R}{X^*,Y}_{\Q,\rho} := \set{Tf:\ f \in
  \Q_{m,n;\rho}}$. If $\epsilon > 0$, we write
$\Pc{R}{X^*,Y}_{\Q,\epsilon}$ and $\Q_{m,n;\epsilon}$ instead of
$\Pc{R}{X^*,Y}_{\Q,(\epsilon,\mdots, \epsilon)}$ and
$\Q_{m,n;(\epsilon,\mdots, \epsilon)}$.  Next, we
put $$\Pc{R}{X^*,Y}_{\Q} := \bigcup_{\rho \in (0, \infty)^{m+n}}
\Pc{R}{X^*,Y}_{\Q,\rho}.$$ For $n=0$ we just write
$\Pc{R}{X^*}_{\Q,\rho}$ instead of $\Pc{R}{X^*,Y}_{\Q,\rho}$.

The properties described Sections \ref{division}, \ref{weierstrass}
and \ref{blowups} of the algebras $\Q_m(U)$ are easily seen to imply
corresponding properties of the algebras $\Q_{m,n;\rho}$ and
$\Pc{R}{X^*,Y}_{\Q,\rho}$.  Due to Proposition \ref{Q-algebra}, we
need no longer formally distinguish between $f \in \Q_{m,n;\rho}$ and
$Tf \in \Pc{R}{X^*,Y}_{\Q,\rho}$; in particular, the notations in
Sections 7, 8 and 9 of \cite{vdd-spe:genpower} make sense in our
setting. 

\begin{df} 
  \label{basicset} 
  A set $A \subseteq I_{m,n;\rho}$ is called a \textbf{basic
    $\Q_{m,n;\rho}$-set} if there are $f,g_1,\mdots, g_k \in
  \Q_{m,n;\rho}$ such that $$A = \set{z \in I_{m,n;\rho} :\ f(z) = 0,\
    g_1(z) > 0,\mdots,\ g_k(z) > 0}.$$ A \textbf{$\Q_{m,n;\rho}$-set}
  is a finite union of basic $\Q_{m,n;\rho}$-sets.  Note that the
  $\Q_{m,n;\rho}$-sets form a boolean algebra of subsets of
  $I_{m,n;\rho}$.
\end{df}

Given a point $a = (a_1,\mdots, a_{m+n}) \in \RR^{m+n}$ and a choice
of signs $\sigma \in \{-1, 1\}^m$, we let $h_{a,\sigma} :\ \RR^{m+n}
\longrightarrow \RR^{m+n}$ be the bijection given by $$h_{a,\sigma}(z)
:= \left(a_1 + \sigma_1 z_1, \mdots, a_{m+n} + z_{m+n}\right). $$ Note
that the maps $h_{a,\sigma}$ (with $a \in \RR^{m+n}$ and $\sigma \in
\{-1, 1\}^m$) form a group of permutations of $\RR^{m+n}$.

\begin{df} 
  \label{gensemianalytic} 
  A set $X \subseteq \RR^{m+n}$ is \textbf{$\Q_{m,n}$-semianalytic at
    $a \in \RR^{m+n}$} if there is an $\epsilon > 0$ such that for
  each $\sigma \in \{-1, 1\}^m$ there is a $\Q_{m,n;\epsilon}$-set
  $A_{\sigma} \subseteq I_{m,n;\epsilon}$ with $$X \cap
  h_{a,\sigma}(I_{m,n;\epsilon}) = h_{a,\sigma}(A_{\sigma}).$$ A set
  $X \subseteq \RR^{m+n}$ is \textbf{$\Q_{m,n}$-semianalytic} if it is
  $\Q_{m,n}$-semianalytic at every point $a \in \RR^{m+n}$.  For
  convenience, if $X \subseteq \RR^m$ is $\Q_{m,0}$-semianalytic we
  also simply say that $X$ is $\Q_m$-semianalytic.
\end{df}

\begin{nrmk} 
  \label{gensemrmk}
  \begin{enumerate}
  \item If $X, Y \subseteq \RR^{m+n}$ are $\Q_{m,n}$-semianalytic at
    $a$, then so are $X \cup Y$, $X \cap Y$ and $X \setminus Y$.
  \item Let $X \subseteq \RR^{m+n}$ be $\Q_{m,n}$-semianalytic, $a \in
    \RR^{m+n}$ and $\sigma \in \{-1, 1\}^m$. Then the set
    $h_{a,\sigma}(X)$ is $\Q_{m,n}$-semianalytic. Moreover by Lemma
    \ref{Q-addition}(3), for each $\lambda \in (0, \infty)^{m+n}$ the
    set $E_{\lambda}(X)$ is $\Q_{m,n}$-semianalytic, where
    $E_{\lambda} : \RR^{m+n} \longrightarrow \RR^{m+n}$ is defined by
    $E_{\lambda}(z) := (\lambda_1 z_1, \dots, \lambda_{m+n} z_{m+n})$.
  \item If $X \subseteq \RR^n$ is semianalytic, then $X$ is
    $\Q_{0,n}$-semianalytic.
  \end{enumerate}
\end{nrmk}

Below we write $0$ for the point $(0,\mdots, 0) \in \RR^{m+n}$.  The
following lemma is now proved just as in \cite[Section
7]{vdd-spe:genpower} with obvious changes: ``$\R_{\mdots}$-set'' is
replaced by ``$\Q_{\dots}$-set'', $\Pc{R}{X^*,Y}$ by
$\Pc{R}{X^*,Y}_{\Q}$ and the algebras $\R_{m,n,\dots}$ by
$\Q_{m,n;\dots}$.  Also, the results from Sections 4, 5 and 6 there
need to be replaced by the corresponding results of Sections
\ref{division}, \ref{weierstrass} and \ref{blowups} here.  (For
example, we use Proposition \ref{Q-basics}(2) here in place of
Corollary 6.7 there; the other replacements are more straightforward.) 

\begin{lemma} \label{gensembasics}
  \begin{enumerate}
  \item Let $A \subseteq \RR^{m+n}$ be $\Q_{m,n}$-semianalytic at $0$
    and let $\sigma$ be a permutation of $\set{1,\mdots, m}$. Then
    $\sigma(A)$ is $\Q_{m,n}$-semianalytic at $0$.
  \item If $n \geq 1$, then each $\Q_{m,n}$-semianalytic subset of
    $\,\RR^{m+n}$ is also $\Q_{m+1,n-1}$-semianalytic.
  \item Every $\Q_{m,n;\rho}$-set $A \subseteq I_{m,n;\rho}$ is
    $\Q_{m,n}$-semianalytic. \qed
  \end{enumerate}
\end{lemma}

Note that Remark \ref{gensemrmk}(3) and Lemma \ref{gensembasics}(2)
imply in particular that every semianalytic subset of $\RR^{m+n}$ is
$\Q_{m,n}$-semianalytic.  Also, since every $f \in \R_{m,n,\rho}$
extends to a holomorphic function $g:B_\frL(\rho') \into \CC$ for some
$\rho'>\rho$, we see that $\R_{m,n,\rho}^\omega \subseteq
\Q_{m,n;\rho}$, where $\R_{m,n,\rho}^\omega$ consists of all $f \in
\R_{m,n,\rho}$ with natural support.  In particular, every
$\R_{m,n}^\omega$-semianalytic subset of $\RR^{m+n}$ is
$\Q_{m,n}$-semianalytic.

We now consider the system $\Lambda = (\Lambda_p)_{p \in \NN}$, where
$$\Lambda_p := \set{A \subseteq I^p:\ A \text{ is }
  \Q_{p} \text{-semianalytic}}.$$ Note that if $A \subseteq I^p$ is
$\Q_{m,n}$-semianalytic with $m+n=p$, then $A$ is also
$\Q_{p}$-semianalytic by Lemma \ref{gensembasics}(2), so $A \in
\Lambda_p$.  A set $A \subseteq \RR^m$ is called a
\textbf{$\Lambda$-set} if $A \in \Lambda_n$, and $B \subseteq \RR^m$
is called a \textbf{sub-$\Lambda$-set} if there exist $n \in \NN$ and
a $\Lambda$-set $A \subseteq \RR^{m+n}$ such that $B = \Pi_m(A)$.

\begin{prop} \label{mc+om1}
Let $A \subseteq [-1,1]^m$ be a sub-$\Lambda$-set.  Then $[-1,1]^m
\setminus A$ is also a sub-$\Lambda$-set.
\end{prop}

\begin{proof}[Sketch of proof]
  By Theorem 2.7 of \cite{vdd-spe:genpower}, we need to establish
  Axioms (I)-(IV) listed in \cite[Section 2]{vdd-spe:genpower}; the
  first three are straightforward.  For Axiom (IV), the proof proceeds
  almost literally as for \cite[Corollary 8.15]{vdd-spe:genpower},
  with the obvious changes indicated earlier, as well as the
  following: $\Pc{R}{X^*, Y}_{\dots}$ there is replaced by
  $\Pc{R}{X^*,Y}_{\Q,\dots}$ here, and the facts of Section
  \ref{o-minimal} here are used in place of the corresponding facts
  from Section 7 there.  Moreover, note that Lemma 6.1 there goes
  through unchanged here.
\end{proof}

Recall that by the remarks after Definition \ref{real_Q}, $$\Rq =
\left(\RR, <, 0, 1, +, -, \mdot, \left(\tilde{f}:\ f \in
    \Q_{m,0;1}\right)\right).$$ It is clear from Remark
\ref{gensemrmk}(2) that every \textit{bounded} $\ \Q_p$-semianalytic
set, for $p \in \NN$, is quantifier-free definable in $\Rq$.  We are
now ready to prove Theorem A.

\begin{thm}
The expansion $\Rq$ is model complete, o-minimal and admits analytic
cell decomposition.
\end{thm}

\begin{proof}
The theorem follows from the previous remark in view of Propositions
\ref{mc+om1} above and \cite[Corollary 2.9]{vdd-spe:genpower}.  For
analytic cell decomposition, we proceed exactly as in the proof of
Corollary 6.10 of \cite{vdd-spe:multisummable}.
\end{proof}

For Theorem B, we proceed as in Section 9 of \cite{vdd-spe:genpower},
with the following changes: we do not need 9.3 there, and use
Corollaries \ref{bad_substitution_1} and \ref{bad_substitution_2} here
in place of Lemma 9.4 there.  Moreover, we do not need 9.7 and Lemma
9.8 there, and we give a much simpler proof of Lemma 9.9 there:

\begin{lemma}
  \label{comp_inverse}
  Let $0 < f \in \Pc{R}{T^*}_\Q$ with $f(0) = 0$.  Then there exists
  $g \in \Pc{R}{T^*}_\Q$ such that $g > 0$, $g(0) = 0$ and $f(g(T)) =
  T$. 
\end{lemma}

\begin{proof}
  By the hypotheses, there are $\lambda,\alpha > 0$ and $h \in
  \Pc{R}{T^*}_\Q$ such that $f(t) = t^\alpha(\lambda + h(t))$ and
  $h(0) = 0$. We put $\rho := 1/\alpha$ and let $F,H \in
  \Pc{R}{T^*,Y^*}_\Q$ be defined by $F(T,Y) := f(Y)$ and $H(T,Y):=
  h(Y)$.  By Proposition \ref{blowup_homom}, the functions
  $\r^{\rho,\lambda} F(t,y)$ and $\r^{\rho,\lambda} H(t,y)$ belong to
  $\Pc{R}{T^*,Y}_\Q$.  We define
  \begin{equation*}
    \phi(t,y) := (\lambda+y)^\alpha \left(\lambda + \r^{\rho,\lambda}
      H(t,y)\right); 
  \end{equation*}
  then $\r^{\rho,\lambda} F(t,y) = t \cdot \phi(t,y)$, so $\phi \in
  \Pc{R}{T^*,Y}_\Q$ by Proposition \ref{Q-basics}(1).  Moreover, we
  have $\phi(0,0) = \lambda^{\alpha+1}$ and $\frac {\partial
    \phi}{\partial y}(0,0) = \alpha\lambda^\alpha > 0$; hence by the
  implicit function theorem, there is a $\psi \in \Pc{R}{T^*}_\Q$ such
  that $\phi(t,\psi(t)) = 1$.  Therefore, we have $\r^{\rho,\lambda} F
  (t,\psi(t)) = t$, so we take $g(t) := t^{1/\alpha} (\lambda +
  \psi(t))$.
\end{proof}

Using this lemma in place of Lemma 9.9 in \cite{vdd-spe:genpower}, we
finish the proof of Theorem B as it is done there, and we obtain
corresponding corollaries for $1$-dimensional sets definable in
$\Rq$.

\section{Example of a definable family of transition maps}
\label{family}

Let $\xi$ be an analytic vector field in $\RR^2$.  A
\textbf{polycycle} $\Gamma$ of $\xi$ is a cyclically ordered finite
set of singular points $p_0 = p_k, p_{1},\dots, p_{k}, p_{k+1} = p_1$
(with possible repetitions), called \textbf{vertices}, and
trajectories $\gamma_{1},\dots, \gamma_{k}$, called
\textbf{separatrices}, connecting the vertices in the order following
the flow of $\xi$, as in Figure 2.  We assume here that each $p_i$ is
a non-resonant hyperbolic singularity of $\xi$.  For each $i$, we fix
two segments $\Lambda_i^-$ and $\Lambda_i^+$ transverse to $\xi$ and
intersecting the separatrices $\gamma_{i-1}$ and $\gamma_i$,
respectively, close to $p_i$.

\begin{figure}[htbp]
  \label{cap:poly}
  \begin{center}
    \input{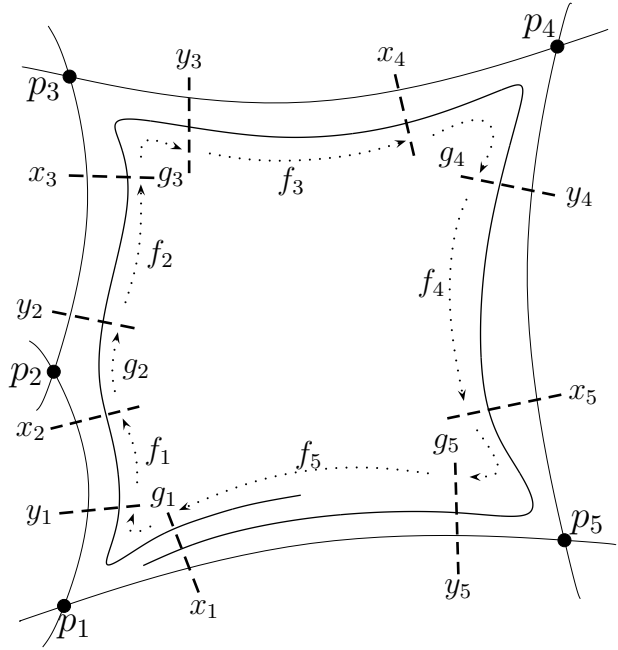}
  \end{center}
  \caption{The polycycle $\Gamma$ of $\xi = \xi_0$}
\end{figure}

For each $i$, we fix analytic charts $x_i:(-1,1) \into \Lambda_i^-$
and $y_i:(-1,1) \into \Lambda_i^+$ such that $x_i(0)$ and $y_i(0)$ are
the points of intersection of $\Lambda_i^-$ with $\gamma_{i-1}$ and of
$\Lambda_i^+$ with $\gamma_i$, respectively, and such that $x_i(t)$
and $y_i(t)$ lie inside the region circumscribed by $\Gamma$ for all
$t \in (0,1)$.  We denote by $g_i:(0,1) \into (0,1)$ the corresponding
transition map in the coordinates $x_i$ and $y_i$; we extend $g_i$ to
all of $(-1,1)$ by putting $g_i(t) := 0$ for $t \in (-1,0]$.  After an
analytic change of coordinates if necessary, it follows from the
general theory of analytic differential equations that there are
analytic functions $f_i:(-1,1) \into (-1,1)$, for $i=1, \dots, k$,
representing the flow of $\xi$ from $\Lambda_{i-1}^+$ to $\Lambda_i^-$
in the charts $y_{i-1}$ and $x_i$.  In fact, these functions $f_i$ are
\textit{restricted analytic}, that is, they extend analytically to a
neighbourhood of $[-1,1]$; in particular, they are definable in $\Rq$.
The restriction of $P := f_k \circ g_k \circ f_{k-1} \circ \cdots
\circ f_1 \circ g_1$ to $(0,1)$ represents the Poincar\'e first return
map of $\xi$ at $p_1$ in the chart $x_1$.

By the corollary and the explanations in the introduction each $g_i$,
and hence $P$, is definable in $\Rq$.  Our goal in this section is to
show that for certain analytic unfoldings $\xi_\mu$ of the vector
field $\xi$, the corresponding Poincar\'e return map \textit{with
  parameter} $\mu$ is also definable in $\Rq$.

More precisely, we let $\xi_{\mu}$ be an analytic unfolding of $\xi$,
with $\mu \in \RR^p$ and $\xi_0 = \xi$, defined in a neighborhood $U$
of $\Gamma$ containing each $\Lambda_i^-$ and $\Lambda_i^+$, with the
same singular points inside $U$ and \textit{with the same linear part
  at each of these singular points as $\xi$}.  We assume that the
unfolding is small, in the sense that for each $\mu \in \RR^p$ and
each $i \in \{1, \dots, k\}$, both $\Lambda_i^-$ and $\Lambda_i^+$
remain transverse to $\xi_\mu$, the transition map of $\xi_\mu$ at
$p_i$ is given by a function $g_{\mu,i}:(0,1) \into (0,1)$ in the
charts $x_i$ and $y_i$ (with the latter as above, independent of
$\mu$), and there are analytic functions $f_{\mu,i}:(-1,1) \into
(-1,1)$ representing the flow of $\xi_\mu$ from $\Lambda_{i-1}^+$ to
$\Lambda_i^-$ in the charts $y_{i-1}$ and $x_i$.  We extend each
$g_{\mu,i}$ to $(-1,1)$ by putting $g_{\mu,i}(t) := 0$ if $t \in
(-1,0]$.  Then the restriction of $P_\mu:= f_{\mu,k} \circ g_{\mu,k}
\circ f_{\mu,k-1} \circ \cdots \circ f_{\mu,1} \circ g_{\mu,1}$ to
$(0,1)$ represents the Poincar\'e first return map of $\xi_\mu$ at
$p_1$ in the chart $x_1$.

We define $g_i:(-1,1) \times \RR^p \into (-1,1)$, $f_i:(-1,1) \times
\RR^p \into (-1,1)$ and $P:(-1,1) \times \RR^p \into (-1,1)$ by
$g_i(t,\mu):= g_{\mu,i}(t)$, $f_i(t,\mu):= f_{\mu,i}(t)$ and
$P(t,\mu):= P_\mu(t)$.  Since the chosen unfolding $\xi_\mu$ is small,
each $f_i$ is a restricted analytic map and hence definable in $\Rq$.
Moreover, we have the following:

\begin{prop}
  \label{finite_cyclicity}
  Each $g_i$ is definable in the structure $\Rq$.  In particular, $P$
  is definable in $\Rq$, and the number of isolated fixed points of
  $P_\mu$ is uniformly bounded in $\mu$.
\end{prop}

For the proof of this proposition, we assume that $p_1 = 0 \in \RR^2$
and show that $g_1$ is definable in $\Rq$.  First, since the ratio
$\lambda = \lambda_2/\lambda_1$ of the eigenvalues $\lambda_1$ (with
respect to $x_1$) and $\lambda_2$ (with respect to $y_1)$ of $\xi_\mu$
at $0$ is irrational and independent of $\mu$, we may assume, after a
change of coordinates that is analytic in both $(x,y)$ and $\mu$, that
the incoming and outgoing separatrices of $\xi_\mu$ at $0$ are
represented by the $x$-axis and the $y$-axis, respectively, for every
$\mu$.  In this situation, the normalization method in
\cite[pp. 70--73]{dul:cycles} goes through uniformly in the parameter
$\mu$ and yields:

\begin{lemma}
  \label{dulac_reduction}
  Let $N \in \NN$ be positive.  Then there exist analytic functions
  $\phi_N,A_N:\RR^{2+p} \into \RR$ such that $A_N(0,0,0) = 0$, the map
  $\Phi^N:\RR^{2+p} \into \RR^{2+p}$ defined by $(u,v,\mu) =
  \Phi^N(x,y,\mu) := (x,\phi_N(x,y,\mu),\mu)$ is a change of coodinates
  fixing $0$ and for each $\mu \in \RR^p$, the push-forward
  $\Phi^N_*\xi_\mu$ satisfies the equations
  \begin{equation}
  \begin{split}
    \label{normal_form}
    \dot{u} & = u\\
    \dot{v} & =
    v\left(\lambda+u^{N}v^{N}A_N\left(u,v,\mu\right)\right). \qed
  \end{split}
  \end{equation}
\end{lemma}

Second, we fix a segment $\Lambda^-:= (0,\epsilon) \times \{y_0\}$,
parametrized by the $x$-coordinate along $\Lambda^-$, and a segment
$\Lambda^+:= \{x_0\} \times (0,\epsilon)$, parametrized by the
$y$-coordinate along $\Lambda^+$.  We assume that $x_0$, $y_0$ and
$\epsilon$ are small enough such that $\Lambda^-$ and $\Lambda^+$ are
transverse to each $\Phi_*^N \xi_\mu$.  We denote by
$g^N_\mu:(0,\epsilon) \into (0,\epsilon)$ the corresponding transition
map of $\Phi_*^N\xi_\mu$, and we define $g^N:(0,\epsilon) \times
\RR^p \into (0,\epsilon)$ by $g^N(t,\mu):= g^N_\mu(t)$.  In this
situation, the estimates obtained on pp. 24--29 in \cite{ily:dulac} go
through uniformly in $\mu$; in particular, the domain $\Omega$ defined
by inequality $( \overset{*}{_{**}} )$ on p. 29 is independent of
$\mu$.  Therefore, we obtain:

\begin{lemma}
  \label{dulac_parameter}
  Let $\nu>0$.  Then there exist an integer $N = N(\nu) > 0$,
  constants $K = K(\nu) > 0$ and $\epsilon = \epsilon(\nu) > 0$ and a
  quadratic domain $\Omega = \Omega(\nu)$ such that
  \begin{enumerate}
  \item the map $g^N$ is analytic in the variable $\mu$ and admits an
    analytic extension to $\Omega \times \RR^p$;
  \item $\left|g^N(t,\mu) - t^{\lambda}\right| \le
    K|t|^{\nu+\epsilon}$ for all $(t,\mu) \in \Omega \times
    \RR^p$. \qed
  \end{enumerate}
\end{lemma}

Third, it follows from the theory of analytic differential equations
that for each $N \in \NN$, there are analytic functions $h_N^-:(-1,1)
\times \RR^p \into (-\epsilon,\epsilon) \times \RR^p$ and
$h_N^+:(-\epsilon,\epsilon) \times \RR^p \into (-1,1) \times \RR^p$
such that $h_N^-(0,\mu) = h_N^+(0,\mu) = 0$ for all $\mu$ and
$g_1(t,\mu) = h_N^+(g^N(h_N^-(t,\mu),\mu),\mu)$ for all $(t,\mu)$.

We write $\langle 1,\lambda \rangle$ for the additive submonoid of
$\RR$ generated by $1$ and $\lambda$.  By the binomial theorem, there
is for each $\alpha \in \langle 1,\lambda \rangle$ and each $N \in
\NN$ an analytic function $c_{\alpha,N}:\RR^p \into \RR$ such that for
each $\mu \in \RR^p$,
\begin{equation*}
  h_N^+\left( (h_N^-(t,\mu))^\lambda,\mu\right) = \sum_{\alpha \in
    \langle 1,\lambda \rangle} c_{\alpha,N}(\mu) \cdot t^\alpha.
\end{equation*}
On the other hand, given $\nu >0$, it follows from Lemma
\ref{dulac_parameter} for each $\mu \in \RR^p$ and each $N \ge N(\nu)$
that
\begin{equation*}
  g_1(t,\mu) - \sum_{\alpha \le \nu} c_{\alpha,N}(\mu) \cdot
  t^\alpha = o\left(\|t\|^\nu\right) \quad\text{as } \|t\| \to 0;
\end{equation*}
in particular, $c_{\alpha,N} = c_{\alpha,N'}$ whenever $|\alpha| \le
\nu$ and $N,N' \ge N(\nu)$.  Thus, for each $\alpha \in \langle
1,\lambda \rangle$ we put $c_\alpha := c_{\alpha,N(|\alpha|)}$; then
by Lemma \ref{dulac_parameter} again, we have for every $\nu > 0$ and
all $(t,\mu) \in \Omega(N(\nu)) \times \RR^p$ that
\begin{equation}
  \label{uni_estimate}
  \left\|g_1(t,\mu) - \sum_{\alpha \le \nu} c_\alpha(\mu) \cdot
    t^\alpha\right\| \le K(\nu) \cdot \|t\|^{\nu+\epsilon(\nu)}.
\end{equation}
It follows from Corollary \ref{uniform_criterion} that $g_1 \in
\A_1(\Omega(0) \times \RR^p)$.  Finally, given any $\nu > \gamma \ge
0$, we define $(g_1)_\gamma:\Omega(N(\gamma)) \times \RR^p \into \RR$
by 
\begin{equation*}
  (g_1)_\gamma(t,\mu) := t^{-\gamma} \left(g_1(t,\mu) - \sum_{\alpha <
      \gamma} 
    c_\alpha(\mu) t^\alpha\right).
\end{equation*}
Then by \eqref{uni_estimate} again, we have for all $(t,\mu) \in
\Omega(N(\nu)) \times \RR^p$ that
\begin{equation*}
  \left\|(g_1)_\gamma(t,\mu) - \sum_{\gamma \le \alpha \le \nu}
    c_\alpha(\mu) \cdot 
    t^{\alpha-\gamma}\right\| \le K(\nu) \cdot
  \|t\|^{\nu+\epsilon(\nu)-\gamma}. 
\end{equation*}
Hence by Corollary \ref{uniform_criterion}, each $(g_1)_\gamma$
belongs to $\A_1(\Omega(N(\gamma)) \times \RR^p)$, that is, $g_1$
satisfies condition (TE).  It follows that $g_1$ belongs to
$\Q_1(\Omega(0) \times \RR^p)$, which proves Proposition
\ref{finite_cyclicity}.


\begin{thebibliography}{10}
\itemsep -0cm

\bibitem{bri-kno:plane}
{\em E.~Brieskorn and H.~Kn\"orrer}, Plane algebraic curves, Birkh\"auser
  Verlag, 1986.

\bibitem{vdd:vol1}
{\em L.~van~den Dries}, Tame Topology and O-minimal Structures, no.~248 in LMS
  Lecture Note Series, Cambridge University Press, 1998.

\bibitem{vdd-spe:genpower}
{\em L.~van~den Dries and P.~Speissegger}, The real field with convergent
  generalized power series is model complete and o-minimal, Trans. Amer. Math.
  Soc., {\bf 350} (1998), ~4377--4421.

\bibitem{vdd-spe:multisummable}
\leavevmode\vrule height 2pt depth -1.6pt width 23pt, The field of reals with
  multisummable series and the exponential function, Proc. London Math. Soc.
  (3), {\bf 81} (2000), ~513--565.

\bibitem{dul:cycles}
{\em H.~Dulac}, Sur les cycles limites, Bull. Soc. Math. France, {\bf 51}
  (1923), ~45--188.

\bibitem{eca:dulac}
{\em J.~Ecalle}, Introduction aux fonctions analysables et preuve constructive
  de la conjecture de Dulac, Hermann, Paris, 1992.

\bibitem{ily:dulac}
{\em Yu.~S. Ilyashenko}, Finiteness theorems for limit cycles, vol.~94 of
  Translations of Mathematical Monographs, American Mathematical Society, 1991.

\bibitem{ily:history}
\leavevmode\vrule height 2pt depth -1.6pt width 23pt, Centennial history of
  {H}ilbert's sixteenth problem, Bull. Amer. Math. Soc., {\bf 39} (2002),
  ~301--354.

\bibitem{rou:bifurcations}
{\em R.~Roussarie}, Bifurcations of planar vector fields and {H}ilbert's
  sixteenth problem, vol.~164 of Progress in Mathematics, Birkh\"auser Verlag,
  Basel, 1998.

\bibitem{spe:pfaffian}
{\em P.~Speissegger}, The {P}faffian closure of an o-minimal structure, J.
  Reine Angew. Math., {\bf 508} (1999), ~189--211.

\bibitem{vdw:algebra}
{\em B.L. van~der Waerden}, Algebra, Springer Verlag, 1967.

\bibitem{whi-wat:analysis}
{\em E.~Whittaker and G.~Watson}, A Course of Modern Analysis, Cambridge, 1927.

\end{thebibliography}
\end{document}